\documentclass[11pt,a4paper]{article}
\RequirePackage{amsmath}%
\RequirePackage{amsthm}%
\usepackage{amsfonts}%
\usepackage{amssymb}%
\usepackage{graphicx}%
\usepackage{float}%
\usepackage{cite}%

\newtheorem{theorem}{Theorem}

\renewcommand{\Re}{\mathop{\mathrm{Re}}}
\renewcommand{\Im}{\mathop{\mathrm{Im}}}

\renewcommand{\i}{\mathrm{i}}

\renewcommand{\Re}{\mathop{\mathrm{Re}}}
\renewcommand{\Im}{\mathop{\mathrm{Im}}}
\renewcommand{\i}{\mathrm{i}}

\newcommand{\CC}{{\mathbb C}}

\newcommand{\bM}{{\bf M}}
\newcommand{\bN}{{\bf N}}
\newcommand{\bI}{{\bf I}}

\newcommand{\bx}{{\bf x}}

\newcommand{\by}{{\bf y}}

\newcommand{\bu}{{\bf u}}
\newcommand{\bU}{{\bf U}}
\newcommand{\bn}{{\bf n}}
\newcommand{\bT}{{\bf T}}
\newcommand{\bZ}{{\bf Z}}

\begin{document}

\title{A fast numerical method for ideal fluid flow in domains with multiple stirrers}

\author{Mohamed M.S. Nasser$^{\rm a}$ and Christopher C. Green$^{\rm b}$}

\date{}
\maketitle

\vskip-0.8cm %
\centerline{$^{\rm a}$Department of Mathematics, Statistics and Physics, College of Arts and Sciences,} %
\centerline{Qatar University, P.O. Box: 2713, Doha, Qatar}%
\centerline{E-mail: mms.nasser@qu.edu.qa}

\centerline{$^{\rm b}$Department of Mathematics, Queensland University of Technology,} %
\centerline{Brisbane, Queensland, Australia}%
\centerline{E-mail: cc.green@qut.edu.au}

\begin{abstract}
A collection of arbitrarily-shaped solid objects, each moving at a
constant speed, can be used to mix or stir ideal fluid, and can give
rise to interesting flow patterns. Assuming these systems of fluid
stirrers are two-dimensional, the mathematical problem of resolving
the flow field - given a particular distribution of any finite number
of stirrers of specified shape and speed - can be formulated as a
Riemann-Hilbert problem. We show that this Riemann-Hilbert problem can
be solved numerically using a fast and accurate algorithm for any
finite number of stirrers based around a boundary integral equation
with the generalized Neumann kernel. Various systems of fluid stirrers
are considered, and our numerical scheme is shown to handle highly
multiply connected domains (i.e. systems of many fluid stirrers) with
minimal computational expense.
\end{abstract}

\section{Introduction}
\label{sc:int}

Riemann-Hilbert (R-H) problems are a classical topic in pure mathematics
that has gradually become more important in applied mathematics~\cite{Abl03,Mus}. 
A R-H problem, following Muskhelishvili's terminology~\cite{Mus},
requires the construction of a function that is analytic everywhere in
a domain in the complex plane satisfying a prescribed boundary
condition on the boundary. Traditionally, R-H problems
arise in the context of singular integral equations, but as has been demonstrated increasingly in recent years, many important problems arising in the applied physical sciences may be cast in the R-H problem framework (e.g., see~\cite{cro-str,skpf,Gak66,Mus,NAS-ETNA15,nas-lap,nas-cr,wegm}). 

Indeed, the main topic of this paper is motivated by a physical
problem. The process of using rigid objects known as stirrers to mix,
or stir, fluid is of interest both from a purely mathematical
standpoint and also from an industrial applications perspective. In a
vessel containing fluid (sometimes dubbed a ‘batch stirring device’),
several fluid stirrers may be inserted and used to stir the fluid.
Each of these stirrers could consist of different, piecewise smooth,
boundary curves and move with different assigned speeds. For modeling
simplicity, we will assume that the fluid being stirred is ideal and
that the stirrers move with specified constant velocities. Since the
stirrers will be assumed to span the depth of the container vessel,
the induced flow can be assumed to be independent of height, and thus
in order to resolve the flow field, the system may be treated as
two-dimensional. This mathematical problem of
determining the instantaneous fluid flow field generated by a
particular assembly of stirrers can be formulated as a R-H
problem defined over a planar multiply connected domain (the shape of
the finite number of stirrers, in addition to the boundary container,
sets the geometry of the fluid domain), and may be readily solved
numerically.

This R-H problem, whose solution will unveil the
instantaneous fluid flow field for a particular collection of
stirrers, may be solved uniquely by a boundary integral equation
having as its kernel function the so-called generalized Neumann kernel
(a generalization of the well-known Neumann kernel). The solvability
of this particular boundary integral equation has been studied 
in~\cite{wegm,nas-cr}. In particular, a fast
and accurate numerical scheme for solving this boundary integral
equation has been presented in~\cite{NAS-ETNA15} and it is this numerical scheme
which we shall employ in this paper in order to solve the problem of
fluid stirrers. The integral equation is discretized according to the
Nystrom method with the trapezoidal rule to obtain a dense linear 
system~\cite{Atk97,Kre14}; this linear system is then solved 
using a combination of the generalized minimal residual (GMRES) 
method~\cite{sad} and the fast multipole
method (FMM) \cite{gre-gim,gre-rok,rok}. 
Such a numerical scheme has been successful in producing solutions to 
an array of problems in conformal mapping and potential theory over
multiply connected domains (see, e.g. 
\cite{Nas-CMFT09,Nas-SIAM09,Nas-JMAA11,Nas-JMAA13,Nas-SIAM13,nas-lap}); 
it has also recently
been shown to expedite the numerical computation of the Schottky-Klein
prime function~\cite{skpf} which is increasingly being demonstrated as an
asset when solving problems in multiply connected domains. Most
noteworthy is the fact that this numerical scheme is able to
accurately solve the boundary integral equation with minimal
computational expense over domains which are highly multiply connected
(as has already been demonstrated in \cite{NAS-ETNA15,nas-jp}); indeed, this will be
exploited to our advantage herein when we consider systems having many
fluid stirrers.

There are several existing mathematical studies on the planar problem
of fluid stirrers. In Wang~\cite{wan} and Burton, Gratus \& Tucker~\cite{bur},
several explicit results were established using conformal mapping
techniques in the case of two circular stirrers. These results were
subsequently generalized by Crowdy, Surana \& Yick~\cite{cro-sur} who catered for
two arbitrarily-shaped moving stirrers, up to knowledge of a conformal
map from a preimage annular domain. Boyland, Aref \& Stremler~\cite{Boy00}
discuss various topological aspects of the viscous flow induced within a
circular disk by three stirrers (thus allowing for chaotic advection
effects), as has also been done by Finn, Cox \& Bryne~\cite{fin} for both
viscous and inviscid fluids for up to five stirrers. The most general
problem of any finite number of arbitrarily-shaped finite stirrers has
been solved by Crowdy~\cite{cro-str} who has constructed an explicit integral
formula, whose kernel function is expressed in terms of the
Schottky-Klein prime function, for the complex potential of this
system, up to knowledge of the conformal mapping from a preimage
multiply connected circular domain. He has highlighted the problem of
fluid stirrers turns out to be equivalent to a modified Schwarz
problem. In the proceeding sections, we will make connections with
each of the works~\cite{cro-str} and~\cite{fin} by recovering certain results
presented therein.

Our paper has the following structure. In section 2, the
R-H problem of the type we will solve is introduced, and a
description given of the aforementioned numerical scheme used to find
solutions to the associated boundary integral equation with the
generalized Neumann kernel. In section 3, the problem of fluid
stirrers is formulated in its most general form as a R-H
problem which is then solved to reveal the instantaneous flow fields
for various configurations of fluid stirrers whose boundaries are
made-up of piecewise smooth curves; in doing so, we recover as special
cases some existing results in~\cite{cro-str} and~\cite{fin}. 
In section 4, we appeal to
conformal slit mappings to study a range of domains having slit, or
paddle type, stirrers. In this section, we also present a novel way of
overcoming the conformal mapping parameter problem based on the work of Aoyama,
Sakajo \& Tanaka~\cite{Aoy-Sak-Tan13} by establishing elliptical or quasi-elliptical preimage
domains (the reason for which will be discussed). Throughout, we
advocate our numerical scheme from the viewpoint of accurately
resolving the flow field around a high number of stirrers (i.e. when
the domains are highly multiply connected) with relative ease and low
computational cost. Concluding remarks are made in section~5.

\section{The boundary integral equation}
\label{sc:aux}

\subsection{The fluid domain $G$}

\begin{figure}[ht] %
\centerline{
\scalebox{0.55}{\includegraphics{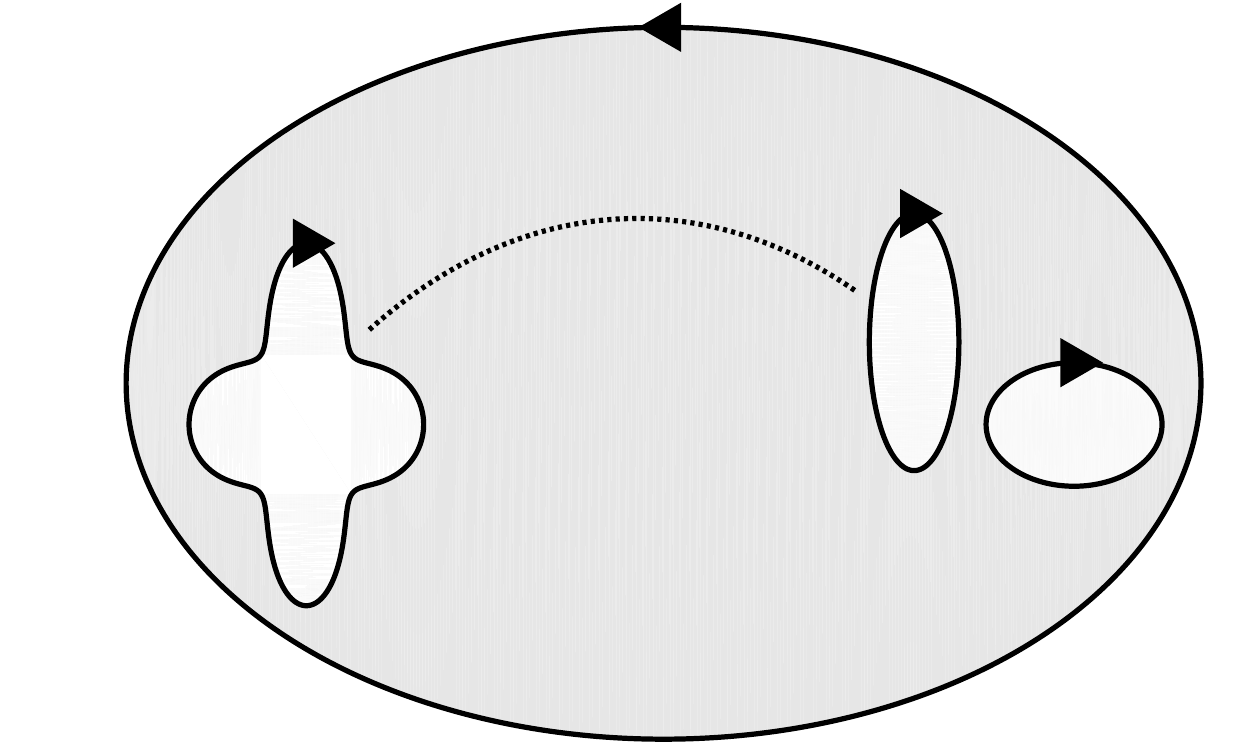}}
\hfill
\scalebox{0.55}{\includegraphics{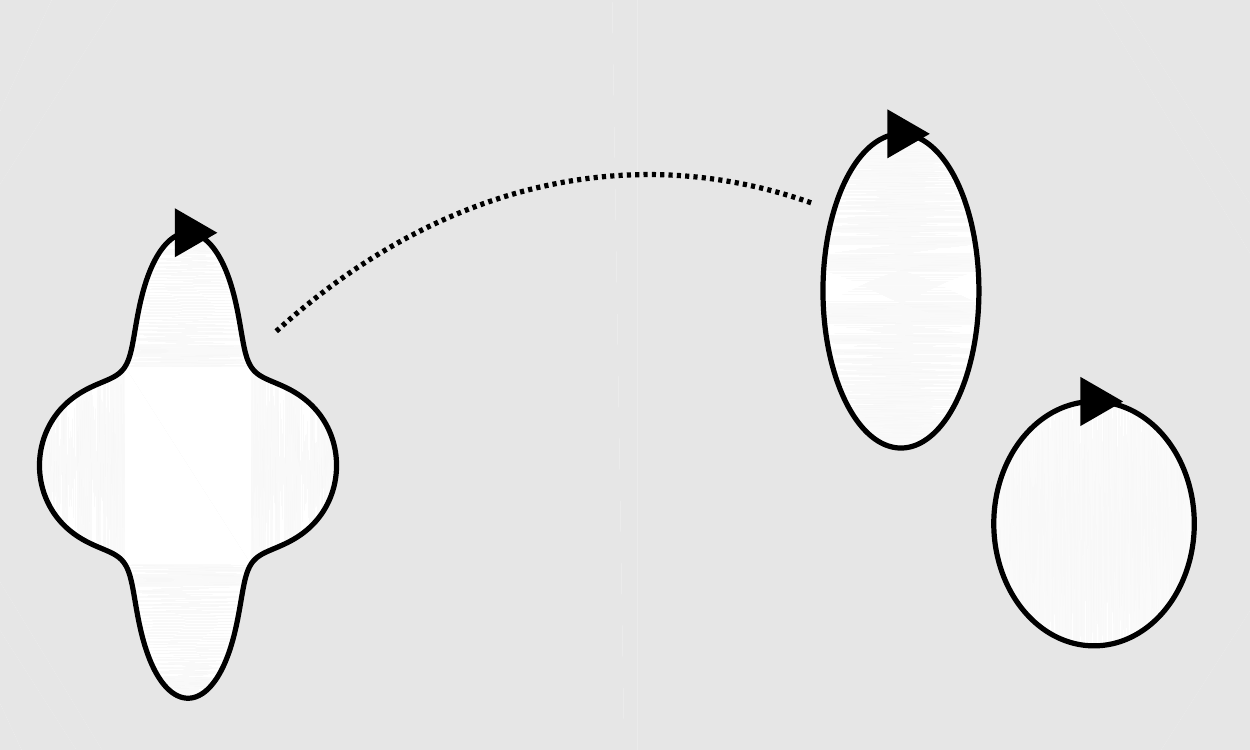}}
}
 \noindent 
 \vskip-4.10cm  \noindent\hspace{3.55cm}  $\Gamma_0$
 \vskip-0.50cm  \noindent\hspace{12.60cm} $\Gamma_1$
 \vskip+0.05cm  \noindent\hspace{8.55cm} $\Gamma_m$
 \vskip-0.55cm  \noindent\hspace{4.90cm}  $\Gamma_2$
 \vskip-0.35cm  \noindent\hspace{1.45cm}  $\Gamma_m$
 \vskip+0.20cm  \noindent\hspace{5.80cm}  $\Gamma_1$
 \vskip-0.25cm  \noindent\hspace{13.65cm} $\Gamma_0$
 \vskip+0.20cm  \noindent\hspace{3.5cm}  $G$  \hspace{7.00cm}  $G$
 \vskip+1.25cm
\caption{ Schematic of a multiply connected domain $G$ of connectivity $m+1$, bounded (left) and unbounded (right).} 
\label{f:B}
\end{figure}

Let $G$ be a multiply connected domain of connectivity $m+1$ ($m\ge0$) in the extended complex plane $\CC\cup\{\infty\}$ having the boundary $\Gamma=\partial G=\cup_{j=0}^{m}\Gamma_j$ where each of the $\Gamma_j$ are closed smooth Jordan curves. The domain $G$ can be either bounded or unbounded. When $G$ is bounded, the curve $\Gamma_0$ encloses the other $m$ curves (see Figure~\ref{f:B}). We assume that $0\in G$ for both bounded and unbounded domains. The orientation of $\Gamma$ will be such that $G$ always lies on the left of $\Gamma$. 

The curve $\Gamma_j$ is parametrized by a $2\pi$-periodic twice continuously differentiable complex function $\eta_j(t)$ with non-vanishing first derivative $\eta'_j(t)\ne 0$, $t\in J_j=[0,2\pi]$, for $j=0,1,\ldots,m$. 
Let $J$ be the disjoint union of $m+1$ intervals $J_0,J_1,\ldots,J_m$ which is defined by 
\begin{equation}\label{e:J}
J = \bigsqcup_{j=0}^{m} J_j=\bigcup_{j=0}^{m}\{(t,j):\;t\in J_j\}.
\end{equation}
The elements of $J$ are ordered pairs $(t,j)$ where $j$ is an auxiliary index indicating which of the intervals the point $t$ lies in. Thus, the parametrization of the whole boundary $\Gamma=\partial G=\Gamma_0\cup\Gamma_1\cup\cdots\cup\Gamma_m$ is defined as the complex function $\eta$ defined on $J$ by
\begin{equation}\label{e:eta-1}
\eta(t,j)=\eta_j(t), \quad t\in J_j,\quad j=0,1,\ldots,m.
\end{equation}
We assume for a given $t$ that the auxiliary index $j$ is known, so we replace the pair $(t,j)$ in the left-hand side of~(\ref{e:eta-1}) by $t$, i.e., for a given point $t\in J$, we always know the interval $J_j$ that contains $t$. The function $\eta$ in~(\ref{e:eta-1}) is thus 
\begin{equation}\label{e:eta}
\eta(t)= \left\{ \begin{array}{l@{\hspace{0.5cm}}l}
\eta_{0}(t),&t\in J_{0}=[0,2\pi],\\
\eta_{1}(t),&t\in J_{1}=[0,2\pi],\\
\hspace{0.3cm}\vdots\\
\eta_m(t),&t\in J_m=[0,2\pi].
\end{array}
\right .
\end{equation}

Let $H$ denote the space of all real functions $\gamma$ defined on $J$, whose restriction $\gamma_j$ to $J_j=[0,2\pi]$ is a
real valued, $2\pi$-periodic and H\"older continuous function for each $j=0,1,\ldots,m$, i.e.,
\[
\gamma(t) = \left\{
\begin{array}{l@{\hspace{0.5cm}}l}
 \gamma_0(t),     & t\in J_0, \\
 \gamma_1(t),     & t\in J_1, \\
  \vdots       & \\
 \gamma_m(t),     & t\in J_m. \\
\end{array}%
\right.
\]
In view of the smoothness of the parametrization $\eta$, a real H\"older continuous function $\hat\gamma$ on $\Gamma$ can be interpreted via $\gamma(t)=\hat\gamma(\eta(t))$, $t\in J$, as a function $\gamma\in H$, and vice versa. So, in this paper, for any given complex or real valued function $\phi$ defined on $\Gamma$, we shall not distinguish between $\phi(t)$ and $\phi(\eta(t))$.

\subsection{The generalized Neumann kernel}

In this paper, we assume the function $A(t)$ is defined for $\eta(t)\in\Gamma$ by
\begin{equation}\label{e:A}
A(t) = \left\{
\begin{array}{l@{\hspace{0.5cm}}l}
 e^{\i\left(\frac{\pi}{2}-\theta(t)\right)}\,(\eta(t)-\alpha),     & {\rm if}\; G\; {\rm is\; bounded}, \\[6pt]
 e^{\i\left(\frac{\pi}{2}-\theta(t)\right)},                    & {\rm if}\; G\; {\rm is\; unbounded}, \\
\end{array}%
\right.
\end{equation}
where $\theta$ is a piecewise constant real valued function defined on $\Gamma$, i.e.,
\begin{equation}\label{e:thet-1}
\theta(t)=\theta_j,\quad t\in J_j,
\end{equation}
where $\theta_j$ are given real constants for $j=0,1,\ldots,m$. 
For simplicity, the piecewise constant function $\theta(t)$ defined on $\Gamma$ by~(\ref{e:thet-1}) will be denoted by
\begin{equation}\label{e:thet}
\theta(t)=(\theta_0,\theta_1,\theta_2,\ldots,\theta_m).
\end{equation}
This notation will also be adopted for any piecewise constant function defined on $\Gamma$.

The generalized Neumann kernel formed with the function $A$ is defined by
\begin{equation}\label{e:N}
 N(s,t) =  \frac{1}{\pi}\Im\left(
 \frac{A(s)}{A(t)}\frac{\eta'(t)}{\eta(t)-\eta(s)}\right).
\end{equation}
We define also a real kernel $M$ by
\begin{equation}\label{e:M}
 M(s,t) =  \frac{1}{\pi}\Re\left(
 \frac{A(s)}{A(t)}\frac{\eta'(t)}{\eta(t)-\eta(s)}\right).
\end{equation}
The kernel $N$ is continuous and the kernel $M$ has a cotangent type singularity.
Hence, the operator $\bN$ defined on $H$ by
\begin{equation}\label{e:opN}
  \bN \mu(s) = \int_J N(s,t) \mu(t) dt, \quad s\in J,
\end{equation}
is a Fredholm integral operator, and the operator $\bM$ defined on $H$ by
\begin{equation}\label{e:opM}
  \bM\mu(s) = \int_J  M(s,t) \mu(t) dt, \quad s\in J,
\end{equation}
is a singular integral operator. For more details, the reader is referred to~\cite{wegm}.

\subsection{The Riemann--Hilbert problem }

Let $\gamma\in H$ be a given function. Following Muskhelishvili~\cite{Mus}, the R-H problem requires determining a function $f(z)$ analytic in $G$, continuous on $G\cup\Gamma$, and $f(\infty)=0$ for unbounded $G$ such that the boundary values of $f$ on $\Gamma$ satisfy
\begin{equation}\label{e:RH-1}
\Re[Af]=\gamma.
\end{equation}
If $A\equiv1$, the R-H problem~(\ref{e:RH-1}) is known as a Schwarz problem or a modified Dirichlet problem~\cite{Gak66,Mus}.

For the function $A$ defined by~(\ref{e:A}), the R-H problem~(\ref{e:RH-1}) is not necessarily solvable for general given function $\gamma\in H$. However, it is always possible to find a unique piecewise constant real $h\in H$ such that the following R-H problem 
\begin{equation}\label{eq:RHp2}
\Re [A f]=\gamma+h
\end{equation}
is solvable~\cite{Nas-CMFT09}. Hence, solving the R-H problem~(\ref{eq:RHp2}) requires determining both the analytic function $f$ as well as the piecewise constant real function $h$. This can be done using a boundary integral equation with the generalized Neumann kernel as in the following theorem established in~\cite{Nas-CMFT09}.

\begin{theorem}\label{thm:rhp}
For a given function $\gamma\in H$, there exists a unique real piecewise constant function of the form
\begin{equation}\label{e:piece}
h(t) = \left\{
\begin{array}{l@{\hspace{0.5cm}}l}
 h_{0},     &t\in J_{0}, \\
 h_{1},     &t\in J_{1}, \\
 \vdots   & \\
 h_m,     &t\in J_m, \\
\end{array}%
\right.
\end{equation}
with real constants $h_{0},h_{1},\ldots,h_m$, such that the R-H problem 
\begin{equation}\label{e:RH-p}
\Re[Af]=\gamma+h
\end{equation}
has a unique solution $f(z)$. The boundary values of the function $f(z)$ are given by 
\begin{equation}\label{e:Af}
Af=\gamma+h+\i\mu
\end{equation}
and the function $h$ is given by
\begin{equation}\label{e:h}
h=[\bM\mu-(\bI-\bN)\gamma]/2
\end{equation}
where $\mu$ is the unique solution of the integral equation
\begin{equation}\label{e:ie}
(\bI-\bN)\mu=-\bM\gamma.
\end{equation}
\end{theorem}

Once the functions $\mu$ and $h$ are found, the boundary values of $f(z)$ follow directly from~(\ref{e:Af}); furthermore, it then follows that the values of $f(z)$ for $z\in G$ can be obtained using the Cauchy integral formula:
\begin{equation}\label{e:f-cau}
f(z)=\frac{1}{2\pi\i}\int_{\Gamma} \frac{\gamma+h+\i\mu}{A} \frac{1}{\eta-z}d\eta.
\end{equation}

\subsection{Numerical solution of the integral equation~(\ref{e:ie})}
\label{sec:num}

In our numerical computations, the boundary integral equation~(\ref{e:ie}) is solved using the MATLAB function $\mathtt{fbie}$ presented in \cite{NAS-ETNA15}. In \verb|fbie|, the integral equation~(\ref{e:ie}) is discretized by the Nystr\"om method with the trapezoidal rule \cite{Atk97,Kre14}. 
Let $n$ be a given even positive integer. Each interval $J_k$ is discretized by the $n$ equidistant nodes
\[
s_{k,p}=(p-1)\frac{2\pi}{n}\in J_k, \quad p=1,2,\ldots,n.
\]
Hence, the total number of nodes in the parameter domain $J$ is $(m+1)n$.
We denote these nodes by $t_i$, $i=1,2,\ldots,(m+1)n$, i.e.,
\begin{equation}\label{e:ti-skp}
t_{(k-1)n+p}=s_{k,p}\in J, \quad k=1,2,\ldots,(m+1)n, \quad p=1,2,\ldots,n.
\end{equation}
Upon specifying a domain with piecewise smooth boundaries, singularity subtraction~\cite{rat} and the trapezoidal rule with a graded mesh~\cite{Kre90} are used.
Hence, we obtain an $(m+1)n \times (m+1)n$ linear algebraic system of the form $(I-B)\bx=\by$. An explicit formula for the elements of the matrix $B$ is given in~\cite{NAS-ETNA15}. This system can be solved iteratively using the MATLAB function $\mathtt{gmres}$. Each step of this method requires one multiplication by the matrix $I-B$. Due to the structure of the integral equation~(\ref{e:ie}), this product is computed efficiently in $O((m+1)n)$ operations using the MATLAB function $\mathtt{zfmm2dpart}$ in the MATLAB toolbox $\mathtt{FMMLIB2D}$ developed by Greengard \& Gimbutas~\cite{gre-gim}. In this way the integral equation~(\ref{e:ie}) is solved in $O((m+1)n\ln n)$ operations. 
In the function $\mathtt{fbie}$, we choose $\mathtt{iprec}=5$ (the tolerance of the FMM 
is $0.5\times 10^{-15}$), $\mathtt{restart}=[\,]$ (the GMRES is used without restart), $\mathtt{gmrestol}=10^{-14}$ (the tolerance of the GMRES method is $10^{-14}$), and $\mathtt{maxit}=100$ (the maximum number of GMRES iterations is $100$). 
For fast numerical evaluation of the Cauchy integral formula~(\ref{e:f-cau}), we use the MATLAB function $\mathtt{fcac}$ in \cite{NAS-ETNA15} which is based on using $\mathtt{zfmm2dpart}$. The method requires $O(\hat n+(m+1)n)$ operations to compute the Cauchy integral formula at $\hat n$ interior points.
The reader is referred to~\cite{NAS-ETNA15} for more details.

\section{Fluid stirrers in domains with piecewise smooth boundaries}
\label{sc:prob}

We consider an incompressible, inviscid, irrotational fluid flow in the domain $G$. The boundary components $\Gamma_j$ are the fluid stirrers whose boundary shapes are specified \textit{a priori}. If $G$ is bounded, then $\Gamma_0$ is the boundary of the fluid vessel. Let $w(z)$ be the complex potential, and hence 
\[
\overline{w'(z)}=u(x,y)+\i v(x,y) 
\]
is the complex velocity of the fluid, and where $\bu=(u,v)$ is its associated velocity field.
On the stirrer $\Gamma_j$ for $j=0,1,\ldots,m$, we have~\cite{cro-str,fin}
\begin{equation}\label{e:u.n}
\bu\cdot\bn_j=\bU_j\cdot\bn_j, 
\end{equation}
where $\bU_j=(U_j,V_j)$ is the specified constant velocity on the stirrer $\Gamma_j$ and $\bn_j$ is the outward-pointing unit normal vector on $\Gamma_j$. For bounded $G$, we define $\bU_0={\bf 0}$; that is, we impose the no-penetration condition~\cite{fin}
\begin{equation}\label{e:u.n-0}
\bu\cdot\bn_0=0,
\end{equation}
which means condition~(\ref{e:u.n}) is satisfied for all $j=0,1,\ldots,m$ for both bounded and unbounded $G$.

Let 
\[
\bT(\eta(t))=\frac{\eta'(t)}{|\eta'(t)|}
\] 
be the unit tangent vector on $\Gamma$. Then,
\[
\bn(\eta(t))=-\i\bT(\eta(t))=-\i\frac{\eta'(t)}{|\eta'(t)|}.
\]
Hence, for $\eta(t)\in\Gamma_j$, the condition~(\ref{e:u.n}) can be written as
\begin{equation}\label{e:w'-cd}
\Re[-\i\eta'(t)w'(\eta(t))]=\Re[-\i\overline{U_j}\,\eta'(t)], \quad t\in J_j. 
\end{equation}
By integrating with respect to the parameter $t$, we obtain~\cite{cro-str}
\begin{equation}\label{e:w-cd}
\Re[-\i w(\eta(t))]=\Re[-\i\overline{U_j}\,\eta(t)]+h_j, \quad t\in J_j,  
\end{equation}
where $h_j$ are real constants of integration. 
The complex potential $w(z)$ will in general be multi-valued.
Let $a_j$ be any point interior to $\Gamma_j$, $j=0,1,\ldots,m$. Then the complex potential $w(z)$ has the form
\begin{equation}\label{e:w-hf}
w(z)=\hat f(z)+\sum_{j=0}^m \frac{\chi_j}{2\pi\i}\log(z-a_j), \quad z\in G\cup\Gamma,
\end{equation}
where $\hat f(z)$ is an analytic function in $G$, and $\chi_j$ is the circulation around the stirrer $\Gamma_j$. If $G$ is bounded, we define $\chi_0=0$. We define a function $f(z)$ for $z\in G\cup\Gamma$ by
\begin{equation}\label{e:f-hf}
f(z)=\frac{\hat f(z)-c}{\i\Pi(z)}
\end{equation}
where
$c=\hat f(0)$ for bounded $G$ and $c=\hat f(\infty)$ for unbounded $G$, and where $\Pi(z)$ is an analytic function defined on $G$ by
\begin{equation}\label{e:Pi}
\Pi(z) = \left\{
\begin{array}{@{}l@{\hspace{0.5cm}}l}
z,& \mbox{if $G$ is bounded},\\
1,         & \mbox{if $G$ is unbounded}. \\
\end{array}%
\right.
\end{equation}
Hence, $f(z)$ is analytic in $G$ (with $f(\infty)=0$ for unbounded $G$) and the complex potential $w(z)$ can be written as
\begin{equation}\label{e:w-f}
w(z)=\i\Pi(z)f(z)+c+\sum_{j=0}^m \frac{\chi_j}{2\pi\i}\log(z-a_j), \quad z\in G\cup\Gamma.
\end{equation}
The constant $c$ has no effect on the velocity field and may be set to zero. It is clear from~(\ref{e:w-f}) that determining the function $f(z)$ is sufficient to fully determine the complex potential $w(z)$, provided the circulations $\chi_j$, $j=0,1,\ldots,m$, are known real numbers. 

Using~(\ref{e:w-cd}), it is straightforward to deduce that the function $f(z)$ is a solution of the following R-H problem  
\begin{equation}\label{e:f-rhp}
\Re[A(t)f(\eta(t))] =\gamma(t)+h(t), \quad t\in J,
\end{equation}
which is of the form in~(\ref{e:RH-p}). Here, the function $A$ is as in~(\ref{e:A}) with $\theta(t)\equiv0$, and
\begin{eqnarray}
\label{e:gam-rhp}
\gamma(t)&=&\Re[-\i\overline{U(t)}\eta(t)]+\sum_{j=0}^m \frac{\chi_j}{2\pi}\log|\eta(t)-a_j|, \\
\label{e:U-rhp}
U(t)&=& (U_0,U_1,\ldots,U_m), \\
\label{e:h-rhp}
h(t)&=& (h_0,h_1,\ldots,h_m).
\end{eqnarray}
The function $f(z)$ will be computed by solving the integral equation~(\ref{e:ie}) as explained in \S\ref{sc:aux}. The complex potential $w(z)$ then follows from~(\ref{e:w-f}).

We shall now use the boundary integral equation~(\ref{e:ie}) to compute the streamline distributions for some bounded and unbounded multiply connected domains whose boundaries consist of various piecewise smooth Jordan curves.

\begin{figure} %
\centerline{
\scalebox{0.6}{\includegraphics{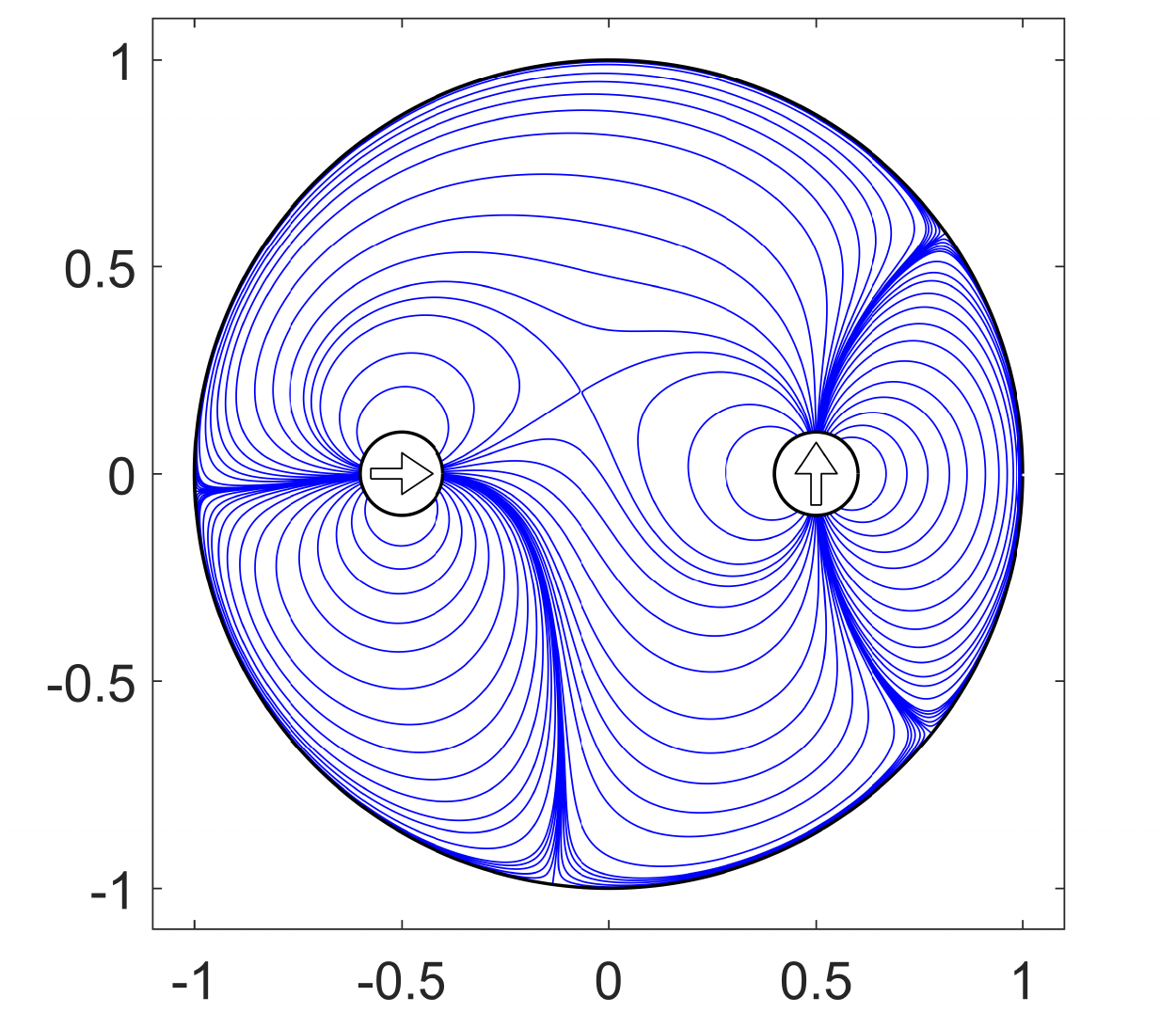}}
\hfill
\scalebox{0.6}{\includegraphics{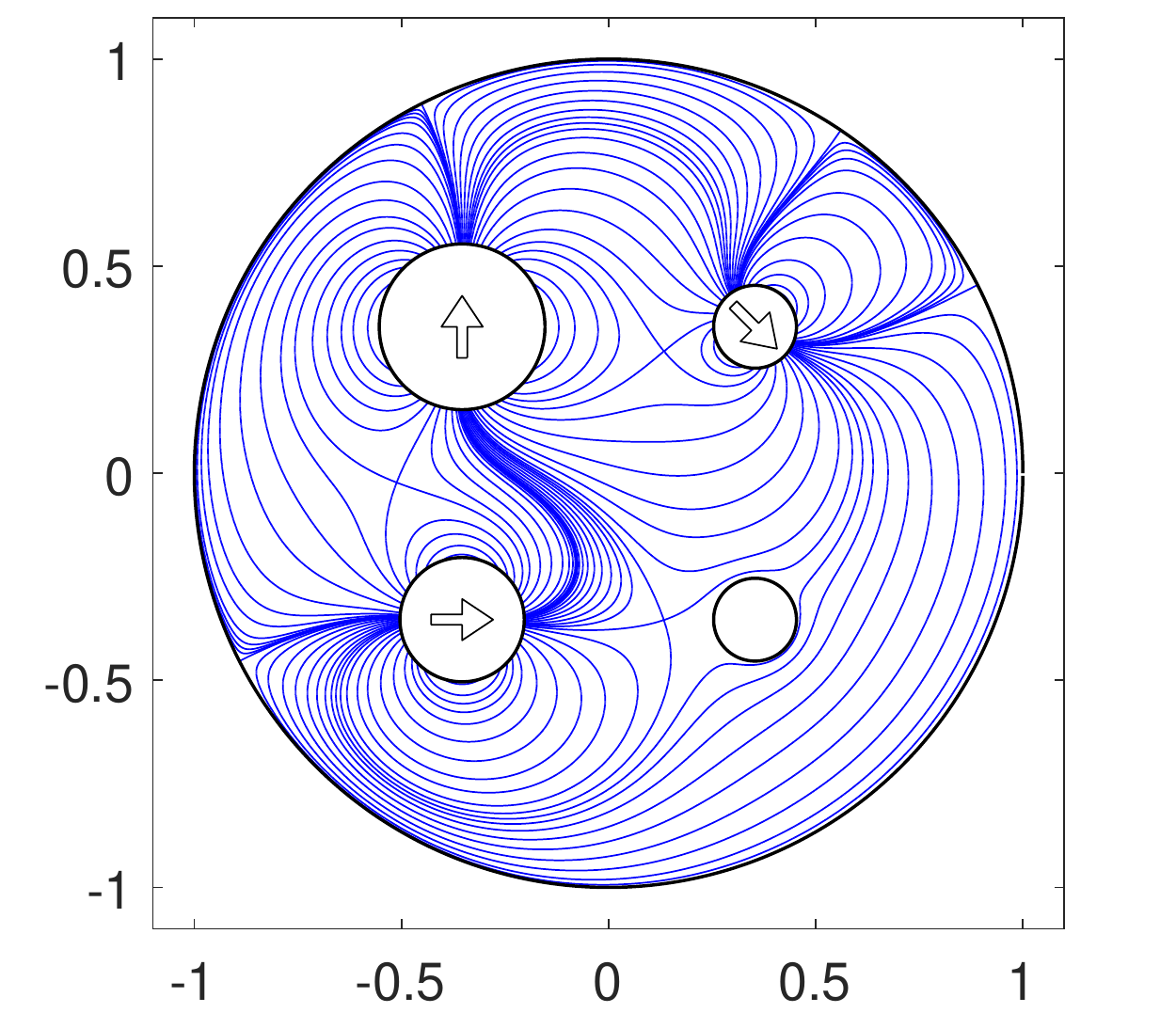}}
}
\caption{Streamlines of the flow generated by two circular disk stirrers (left) and four circular disk stirrers (right) inside the unit disk. The stirrers in the left figure are centered at $-0.5,0.5$, each are of radius $0.1$, and have complex velocities $1,\i$, respectively. The stirrers in the right figure are centered at $0.5e^{3\pi\i/4},0.5e^{\pi\i/4},0.5e^{-\pi\i/4},0.5e^{-3\pi\i/4}$, have radii $0.2,0.1,0.1,0.15$, and have complex velocities $\i,e^{-\pi\i/4},0,1$, respectively.} 
\label{f:disks}
\end{figure}

\begin{figure} %
\centerline{
\scalebox{0.60}{\includegraphics{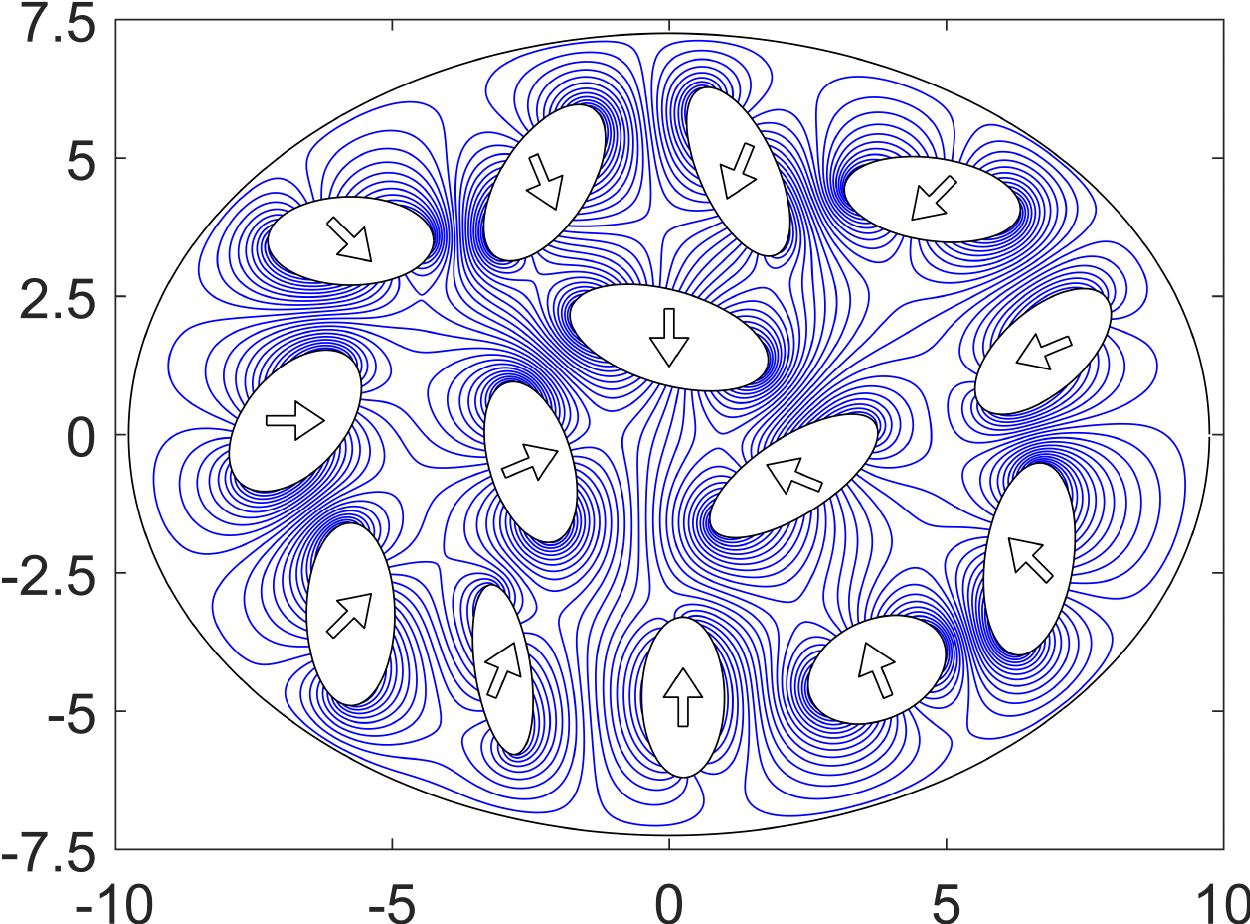}}
\hfill
\scalebox{0.60}{\includegraphics{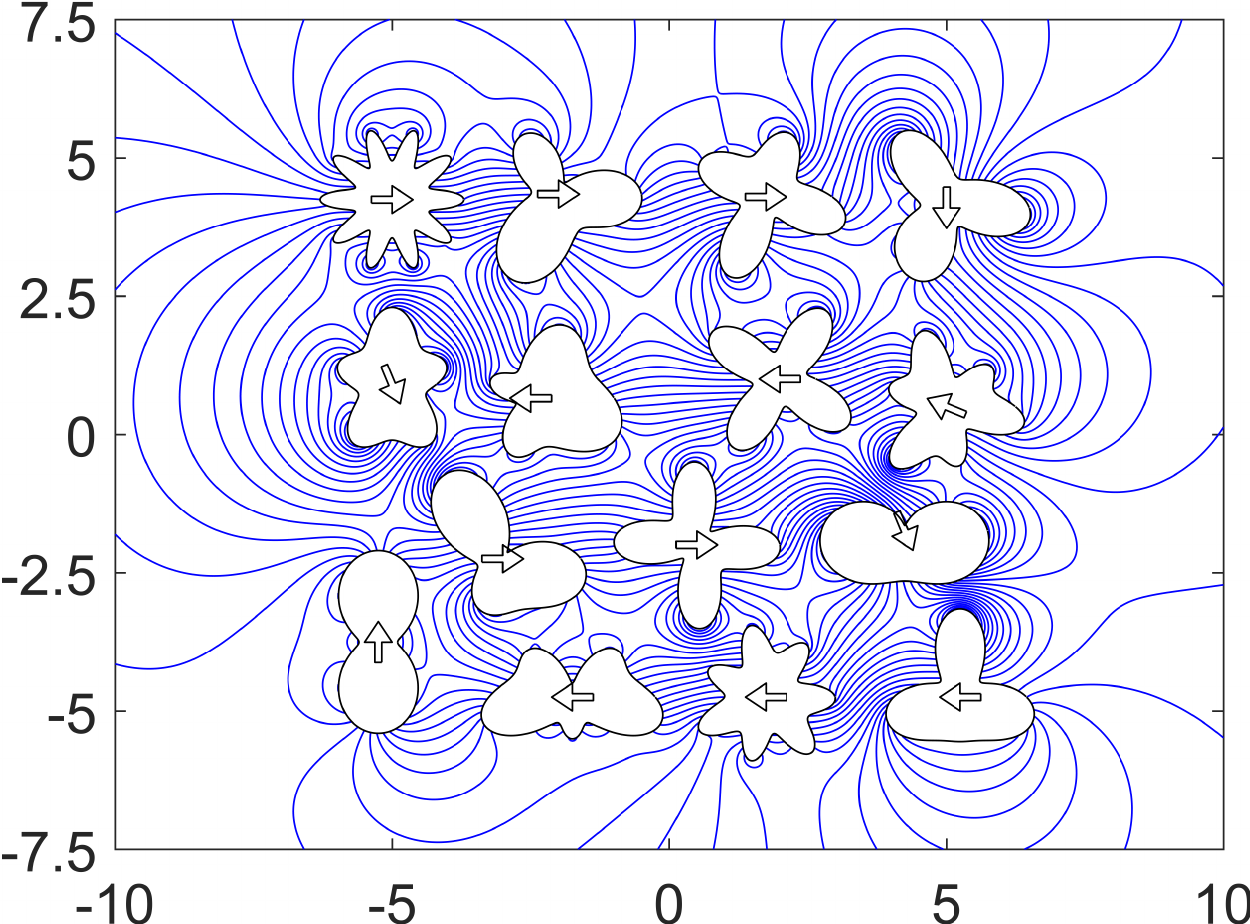}}
}
\caption{Streamlines of the flow generated by fourteen elliptical stirrers inside an ellipse (left) and fifteen stirrers of various size and shape in an unbounded domain (right). All stirrers have random circulations between $-1$ and $+1$.} 
\label{f:bou-unbou-15}
\end{figure}

\begin{figure} %
\centerline{
\scalebox{0.60}{\includegraphics{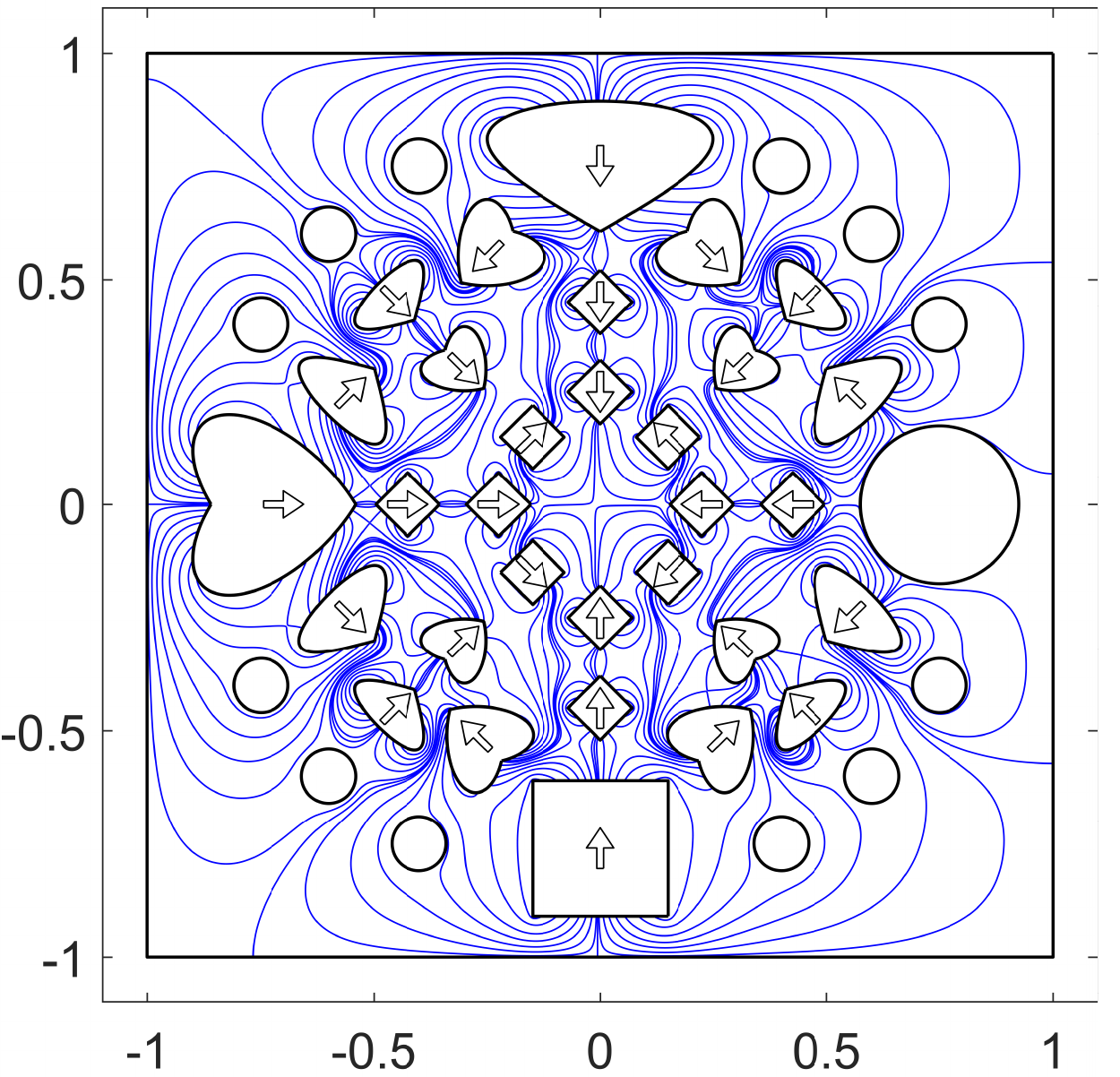}}
\hfill
\scalebox{0.60}{\includegraphics{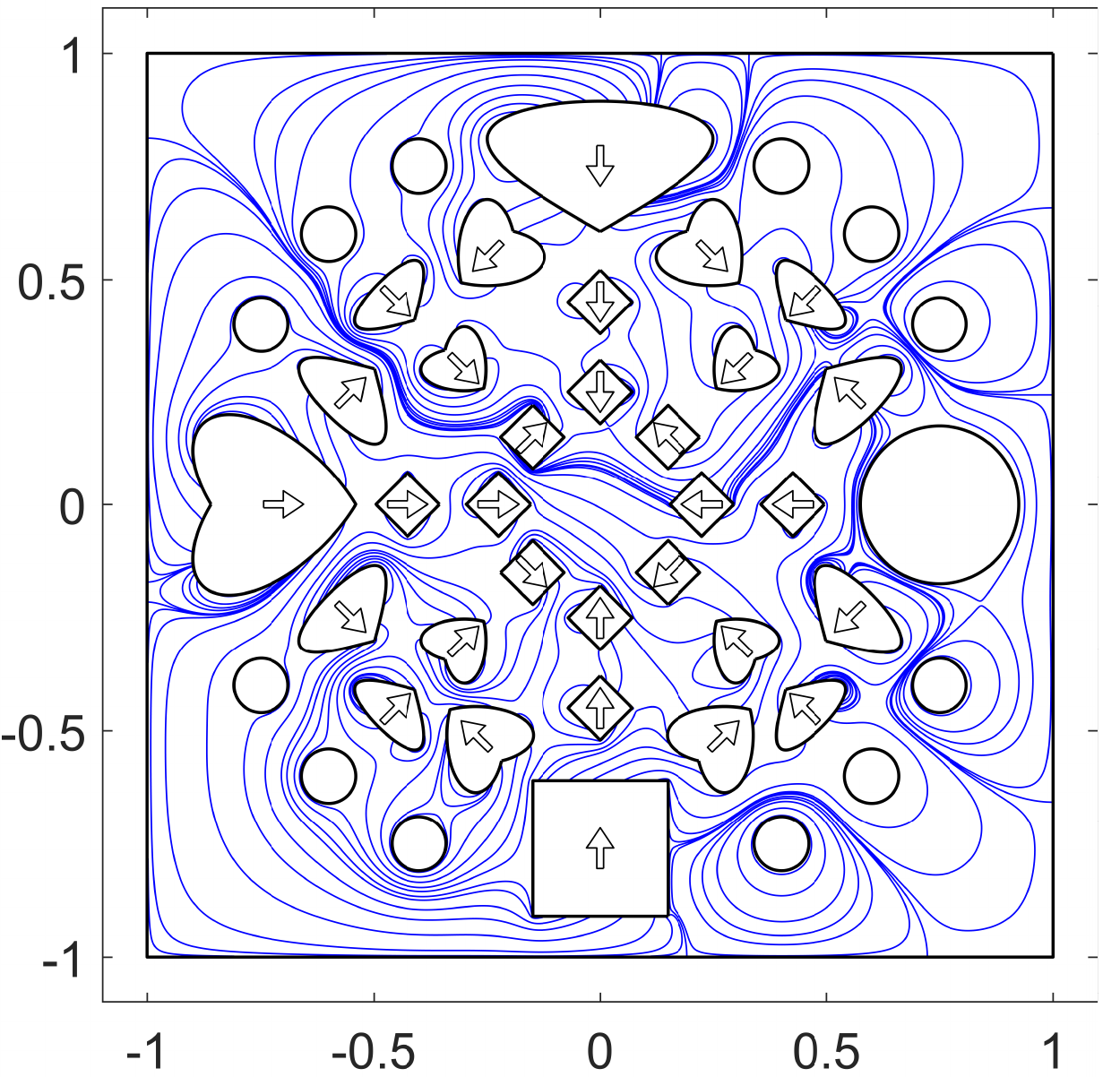}}
}
\caption{Streamlines of the flow generated by forty-four stirrers of various shape and size inside a square. All internal stirrers have zero circulation (left) and random circulations between $-1$ and $+1$ (right).} 
\label{f:Squares-44}
\end{figure}

\begin{figure} %
\centerline{\scalebox{0.95}{\includegraphics{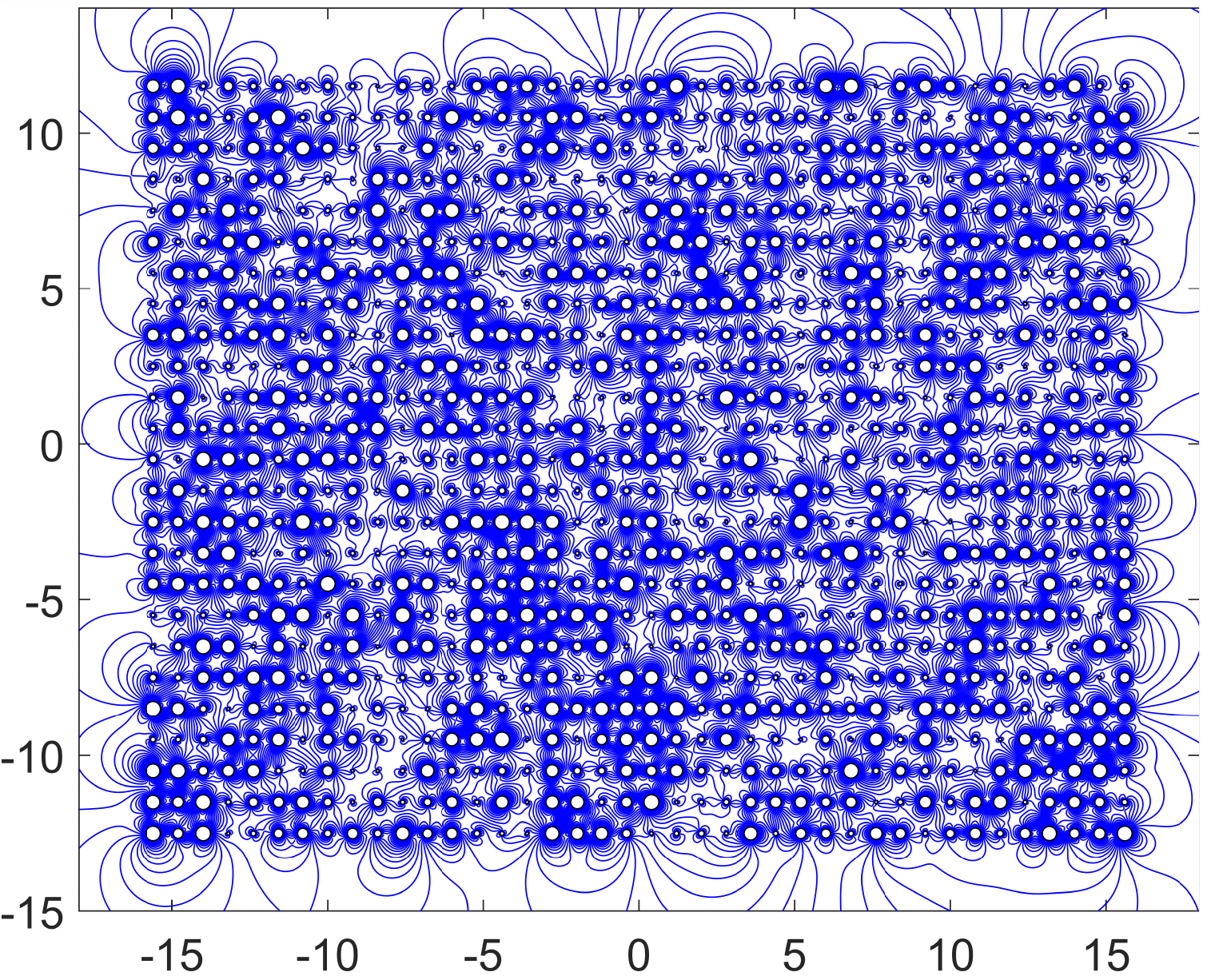}}}
\caption{Streamlines of the flow generated by one thousand circular stirrers in an unbounded domain. The radii of these stirrers are random real numbers between $0.04$ and $0.24$. All stirrers have complex velocities of modulus $1$ in random directions.} 
\label{f:disks-1000}
\end{figure}

\begin{figure} %
\centerline{\scalebox{0.60}{\includegraphics{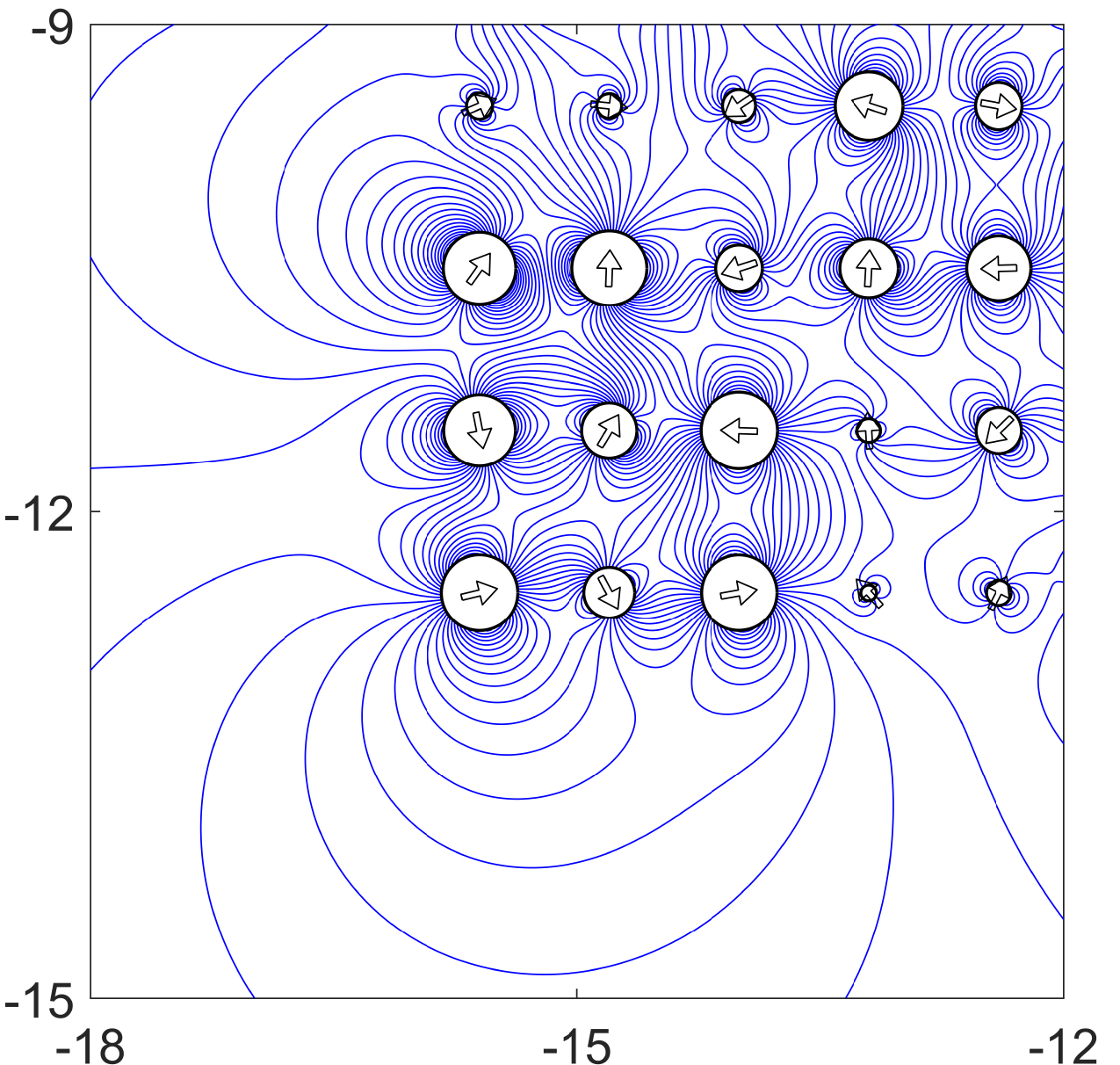}}}
\caption{A magnified sub-section of the lower left-hand domain of the flow domain of~Figure~\ref{f:disks-1000}.} 
\label{f:disks-1000z}
\end{figure}

To check that our numerical scheme recovers some existing results, we computed for Figure~\ref{f:disks} some streamline distributions presented in both Crowdy~\cite{cro-str} and Finn et al.~\cite{fin} for the flow in two bounded multiply connected circular domains (see Figure~\ref{f:disks}). 
We see that ours are in very good qualitative agreement with theirs. 
In Figure~\ref{f:bou-unbou-15}, we show the streamlines of the flow due to fluid stirrers in a bounded and unbounded domain each of connectivity fifteen. 
To demonstrate that our numerical scheme can be used for domains with piecewise smooth Jordan curves, we show in Figure~\ref{f:Squares-44} the streamlines of the flow generated by forty-four stirrers inside a square. Further, to show that our numerical method can also be used effectively for domains with very high connectivity, we show in Figure~\ref{f:disks-1000} the streamlines of the flow due to one thousand circular disk stirrers in an unbounded domain. 
In Figures~\ref{f:disks}--\ref{f:Squares-44}, stirrers with arrows inside them have complex velocity of modulus $1$ in the directions indicated by the arrows. Stirrers with no arrows are stationary.
All stirrers in Figure~\ref{f:disks-1000} have complex velocities of modulus $1$ in arbitrary directions. The directions are indicated by the arrows in Figure~\ref{f:disks-1000z} which is a magnified sub-section of the lower left-hand domain of the flow domain of Figure~\ref{f:disks-1000} showing twenty circular stirrers. 
The fluid stirrers in Figures~\ref{f:disks}, \ref{f:Squares-44} (left) and~\ref{f:disks-1000} have zero circulation around them whilst the fluid stirrers in Figures~\ref{f:bou-unbou-15} and~\ref{f:Squares-44} (right) have been allocated random circulations between $-1$ and $+1$.

The following table shows the computation times (in seconds) for solving the integral equation~(\ref{e:ie}) for the domains in Figures~\ref{f:disks}--\ref{f:disks-1000}. Computation times were measured using the MATLAB \verb|tic toc| command on a standard laptop computer.
\begin{center}
\begin{tabular}{llll}
\hline
Domain & $n$ & Total number of nodes & Time (s) \\
\hline
Figure~\ref{f:disks} (left)         & $1024$ & $3072$    & $0.2$ \\
Figure~\ref{f:disks} (right)        & $1024$ & $5120$    & $0.3$ \\
Figure~\ref{f:bou-unbou-15} (left)  & $1024$ & $15360$   & $1.2$ \\
Figure~\ref{f:bou-unbou-15} (right) & $1024$ & $15360$   & $1.8$ \\
Figure~\ref{f:Squares-44}           & $2048$ & $92160$   & $6.2$ \\
Figure~\ref{f:disks-1000}           & $1024$ & $1024000$ & $36.0$ \\
\hline
\end{tabular}
\end{center}

It is important to note that Crowdy~\cite{cro-str} found explicit formulae for the ideal flow produced by any finite number of fluid stirrers, circular or otherwise. These analytical formulae are written in terms of the Schottky-Klein prime function which should be computed using software based on the numerical schemes of~\cite{skpf}. The fastest method presented in~\cite{skpf} is based on a boundary integral equation of the form in~(\ref{e:ie}) and computing the Schottky-Klein prime function requires solving at least $m+2$ of these integral equations; however, the method presented in this paper requires solving only one such integral equation. Thus, if one were to use the formulae of~\cite{cro-str} to compute the streamlines in Figure~\ref{f:disks-1000}, the total computational cost would be considerably greater compared to that due to our method. Furthermore, the formulae in~\cite{cro-str} can be used for domains with arbitrarily-shaped fluid stirrers (such as the domains in Figures~\ref{f:bou-unbou-15} and~\ref{f:Squares-44}) provided a conformal map from a multiply connected circular domain is known. In this case, in addition to computing the Schottky-Klein prime function, it will also be required to compute the conformal mapping numerically by, for example, the methods presented in~\cite{NAS-CMFT15,weg-fast}, and the computational cost of these methods is much higher than the overall computational cost of the method presented in this section. Indeed, the computational costs of the methods presented in~\cite{NAS-CMFT15} and~\cite{weg-fast} are $O(m^2n+mn\log n)$ and $O(m^2n^2)$, respectively, whereas the computational cost of our method is just $O(mn\log n)$. 

\section{Fluid slit stirrers and conformal mapping}
\label{sc:wall}

The method presented in the previous section can be used to compute the streamlines of the ideal fluid flow generated by stirrers whose boundaries are piecewise smooth Jordan curves. 
However, with the aide of conformal mapping, the method can be extended to include stirrers shaped as slits, which of course are not Jordan curves. 
Let $\Omega$ be a multiply connected domain in the $\zeta$-plane whose boundaries are slits. Suppose that $\Omega$ is the image under a suitable conformal mapping $\zeta=\Phi(z)$ of a multiply connected domain $G$ in the $z$-plane whose boundaries are piecewise smooth Jordan curves (i.e., of the type considered so far). Let $\hat\Gamma_j=\Phi(\Gamma_j)$. Then $\partial\Omega=\hat\Gamma=\cup_{j=0}^{m}\hat\Gamma_j$. 

Suppose that $w(\zeta)$ is the complex potential of the flow in the fluid domain $\Omega$. The function $w(\zeta)$ satisfies on $\hat\Gamma$ the boundary condition
\begin{equation}\label{e:w-cd-zt}
\Re[-\i w(\zeta)]=\Re[-\i\overline{U_j}\,\zeta]+h_j, \quad \zeta\in \hat\Gamma_j, \quad j=0,1,\ldots,m.
\end{equation}
Using $\zeta=\Phi(z)$, we see that $W(z)=w(\Phi(z))$ satisfies on $\Gamma$ the boundary condition
\begin{equation}\label{e:W-cd-zt}
\Re[-\i W(z)]=\Re[-\i\overline{U_j}\,\Phi(z)]+h_j, \quad z\in\Gamma_j, \quad j=0,1,\ldots,m.
\end{equation}
The reader is referred to~\cite{cro-str} for more details.

As in the previous section, the complex potential $W(z)$ can be written as
\begin{equation}\label{e:w-f-W}
W(z)=\i\Pi(z)f(z)+\sum_{j=0}^m \frac{\chi_j}{2\pi\i}\log(z-a_j), \quad z\in G\cup\Gamma.
\end{equation}
where the function $f(z)$ is an analytic function in the domain $G$ with $f(\infty)=0$ for unbounded $G$. The function $f(z)$ is a solution of the RH problem
\begin{equation}\label{e:f-rhp-W}
\Re[A(t)f(\eta(t))] =\gamma(t)+h(t), \quad t\in J,
\end{equation}
where
\begin{eqnarray}
\label{e:gam-rhp-W}
\gamma(t)&=&\Re[-\i\overline{U(t)}\Phi(\eta(t))]+\sum_{j=0}^m \frac{\chi_j}{2\pi}\log|\eta(t)-a_j|, \\
\label{e:U-rhp-W}
U(t)&=& (U_0,U_1,\ldots,U_m), \\
\label{e:h-rhp-W}
h(t)&=& (h_0,h_1,\ldots,h_m),
\end{eqnarray}
and we recall that we set $\chi_0=0$ and $U_0=0$ for bounded $G$.
The function $f(z)$ is found in the same way as before.

By computing the analytic function $f(z)$, we obtain the complex potential $W(z)$. The complex potential $w(\zeta)$ for the slit domain $\Omega$ is then given by $w(\zeta)=W(\Phi^{-1}(\zeta))$. 
We note that in order to compute the streamlines in the slit domain $\Omega$, it is not required to compute the inverse map $\Phi^{-1}$. Instead, we discretize the domain $G$ and the direct mapping $\Phi$ is used to obtain a discretization of the slit domain $\Omega$. Then the values of the function $w(\Phi(z))$ are used to compute the streamlines in $\Omega$.

This method for computing the streamlines for the slit domain $\Omega$ can be summarized as follows:
\begin{itemize}
	\item Compute the preimage domain $G$ and the conformal mapping $\Phi$ from $G$ onto $\Omega$.
	\item Let $\bZ$ be a matrix of points obtained by discretizing the preimage domain $G$ (if $G$ is unbounded, then we discretize only small part of $G$ surrounding the boundaries of $G$). Then $\boldsymbol{\zeta}=\Phi(\bZ)$ are discretizing points of the domain $\Omega$.
	\item Solve the R-H problem~(\ref{e:f-rhp-W}) in the preimage domain $G$ for the analytic function $f(z)$. Hence $W(z)$ is given by~(\ref{e:w-f-W}). We compute the values of the function $W(z)$ at the points $\bZ$. 
	\item Then we compute the values of the function $w(\zeta)$ at the points $\boldsymbol{\zeta}$ through 
	$w(\boldsymbol{\zeta})=W(\Phi^{-1}(\boldsymbol{\zeta}))=W(\bZ)$. Then plot the contour lines of the function $w(\zeta)$.
\end{itemize}

As should be apparent, it is straightforward to compute the streamlines for the slit domain $\Omega$ provided we know the preimage domain $G$ and the conformal mapping $\Phi$ from $G$ onto $\Omega$. However, knowing the preimage domain $G$ and the conformal mapping $\Phi$ from $G$ onto $\Omega$ is not a simple task. One of the main contributions of this paper will be providing a numerical method for computing the preimage domain $G$ and the conformal mapping $\Phi$ from $G$ onto $\Omega$ for a given slit domain $\Omega$. This numerical method will be presented in the remaining of this section.

In this paper, by way of example, we shall consider the following two canonical slit domains $\Omega$:
\begin{itemize}
	\item The entire $\zeta$-plane with $m+1$ finite rectilinear slits.
	\item The upper half-plane with $m$ finite rectilinear slits.
\end{itemize}
The method can be readily extended to cater for other canonical slit domains.

An efficient numerical method for computing the conformal mapping from any given bounded or unbounded multiply connected domain $G$ bounded by Jordan curves onto the above two canonical slit domains and onto more other canonical slit domains has been developed in a series of papers~\cite{Nas-CMFT09,Nas-SIAM09,Nas-JMAA11,Nas-JMAA13}. The method is based on a unified boundary integral equation with the generalized Neumann kernel. In these papers, the domain $G$ is assumed to be known and the integral equation is used to find the conformal mapping as well as the canonical slit domain $\Omega$. 
However, in this paper, we need to compute the streamlines for a given slit domain $\Omega$; that is, we assume that the slit domain $\Omega$ is known. Hence, the preimage domain $G$ will be unknown. Thus, we need to compute the preimage domain $G$ as well as the conformal mapping $\zeta=\Phi(z)$ from $G$ onto $\Omega$.

For the first canonical domain (the entire $\zeta$-plane with $m$ finite rectilinear slits), an iterative numerical method for computing the preimage domain $G$ and the conformal mapping $\zeta=\Phi(z)$ has been suggested in Aoyama, Sakajo \& Tanaka~\cite{Aoy-Sak-Tan13} where the preimage $G$ is assumed to be circular. Since the image domain is elongated (slit domains), numerical crowding effects are problematic. Further, the circles will be close to each other and the iterative method will either be slow to converge or fail to do so altogether. To overcome such difficulties, we shall assume in this paper that the preimage domain $G$ is bounded by ellipses instead of circles. 

\subsection{The entire $\zeta$-plane with $m$ finite rectilinear slits}
\label{sc:rec}

\begin{figure}[ht] %
\centerline{
\scalebox{0.5}{\includegraphics{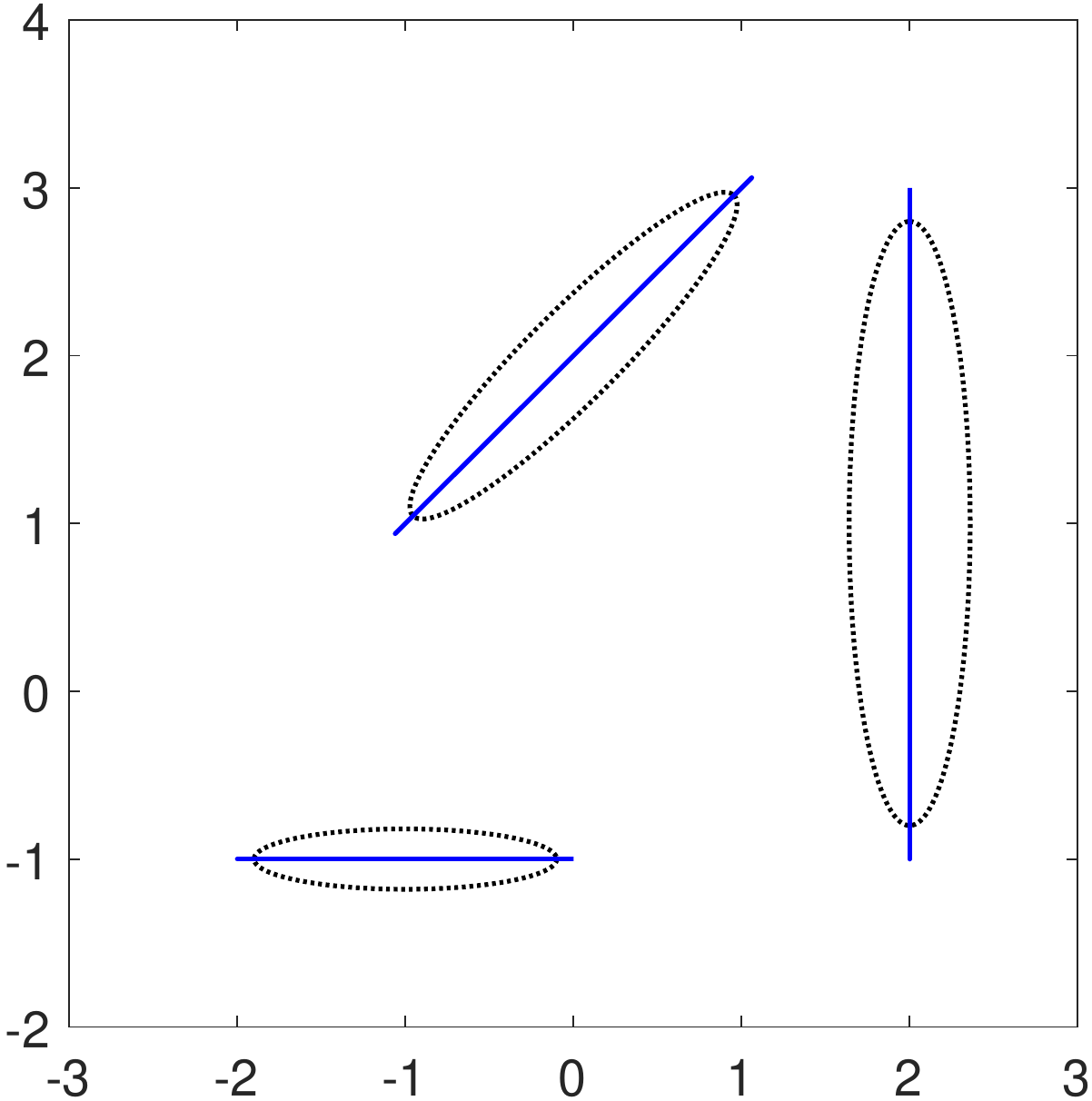}}
\hfill
\scalebox{0.5}{\includegraphics{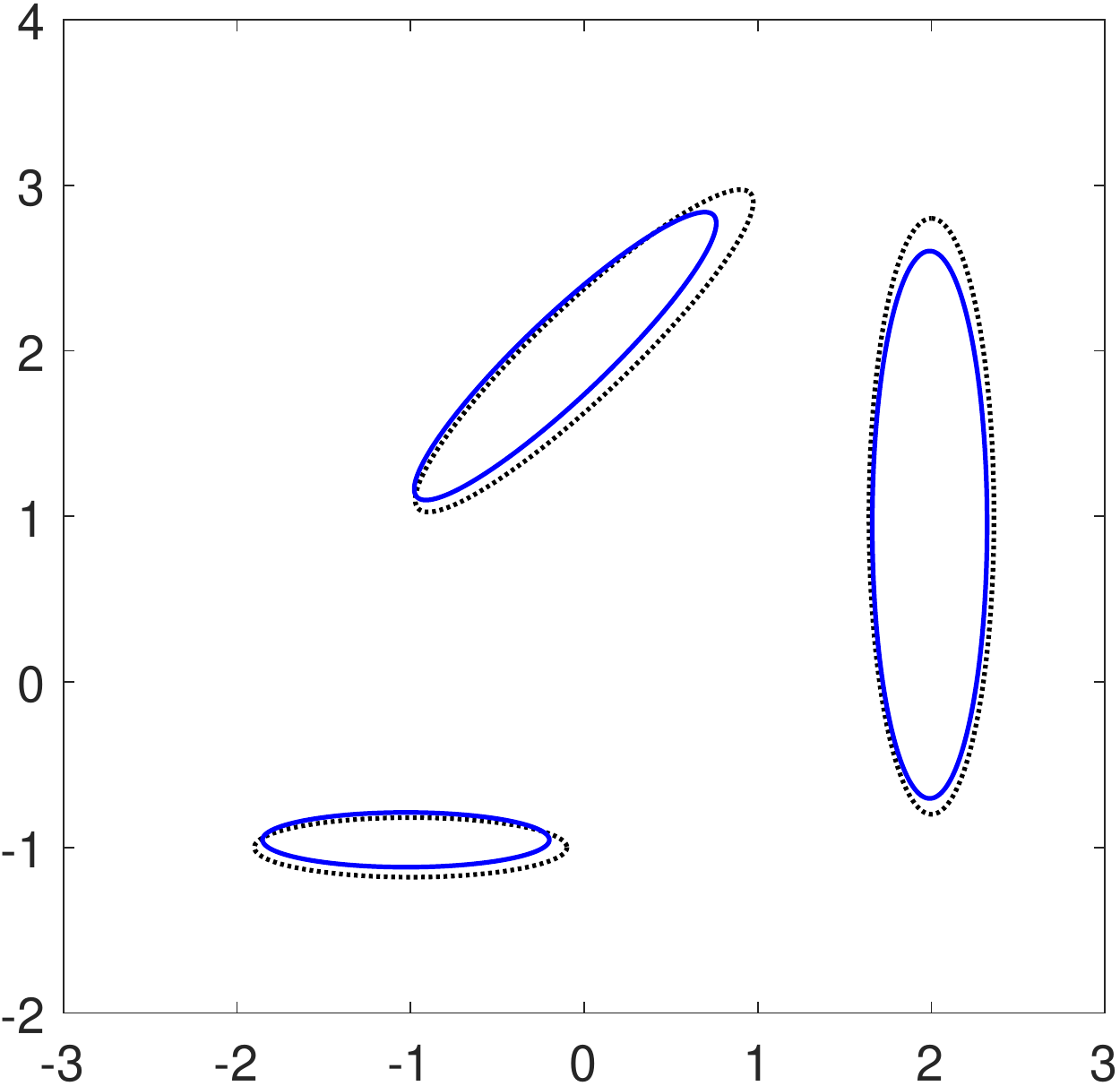}}
}
\caption{A given rectilinear slit domain (solid line) with the initial preimage domain $G^0$ (dotted line) (left); and the initial preimage domain $G^0$ (dotted  line) with the computed preimage domain G (solid line) (right)} 
\label{f:str}
\end{figure}

Let $\Omega$ be the entire $\zeta$-plane with $m+1$ rectilinear slits $L_j$, $j=0,1,\ldots,m$, making angles $\theta_j$ with the positive real line (see Figure~\ref{f:str} (left) for $m=2$). For such canonical domains, we shall assume the preimage domain $G$ is an unbounded multiply connected domain exterior to $m+1$ ellipses. Assuming the boundary $\Gamma=\partial G$ is parametrized as in~(\ref{e:eta}), then the conformal mapping $\zeta=\Phi(z)$ with the normalization
\[
\Phi(\infty)=\infty, \quad \lim_{z\to\infty}(\Phi(z)-z)=0.
\]
can be computed as in the following theorem from~\cite{Nas-JMAA11}.

\begin{theorem}\label{thm:cm-rec-slit}
Let $\theta$ be the piecewise constant function defined on $\Gamma$ by $\theta(t)=(\theta_0,\theta_1,\ldots,\theta_m)$, the function $A$ be defined by~(\ref{e:A}), and the function $\gamma$ be defined by
\begin{equation}
\gamma(t)=\Im\left[e^{-\i\theta(t)}\eta(t)\right], \quad t\in J.
\end{equation}
Let also $\mu$ be the unique solution of the boundary integral equation~(\ref{e:ie}) and the piecewise constant function $h$ be given by~(\ref{e:h}). Then the function $f$ with the boundary values
\begin{equation}\label{eq:f-rec}
f(\eta(t))=(\gamma(t)+h(t)+\i\mu(t))/A(t)
\end{equation}
is analytic in $G$ with $f(\infty)=0$ and the conformal mapping $\Phi$ is 
given by
\begin{equation}
\Phi(z)=z+ f(z), \quad z\in G\cup\Gamma.
\label{eqn:omega-app}
\end{equation}
\end{theorem}

The application of Theorem~\ref{thm:cm-rec-slit} requires that the domain $G$ is known. However, for our case it is $\Omega$ which is known and the domain $G$ is unknown and needs to be determined alongside the conformal mapping  $\zeta=\Phi(z)$ from $G$ onto $\Omega$. An iterative method for computing the preimage domain $G$ and the conformal mapping $\zeta=\Phi(z)$ will be described in this section. 
The iterative method is used to generate a sequence of multiply connected domains $G^0,G^1,G^2,\ldots$ which converge to the required preimage domain $G$. 

Let $\ell_j$ denote the length of the slit $L_j$, let $\zeta_j$ denote its center, and let $\theta_j$ denote the angle of intersection between the line and the positive real axis ($\ell_j$, $\zeta_j$ and $\theta_j$ are given for $j=0,1,\ldots,m$). In the iteration step $k=0,1,2,\ldots$, we assume the boundaries $\Gamma^k_0,\Gamma^k_1,\ldots,\Gamma^k_m$ of the domain $G^k$ are the ellipses parametrized by
\begin{equation}\label{eq:eta-k}
\eta^k_j(t)=z^k_j+0.5e^{\i\theta_j}(a^k_j\cos t-\i b^k_j\sin t), \quad 0\le t\le 2\pi, 
\end{equation}
for $j=0,1,\ldots,m$, where the parameters of these ellipses, i.e. the centers of the ellipses $z^k_j$, the lengths of the major axes $a^k_j$, and the lengths of the minor axes $b^k_j$, will be computed using the following iterative method which is a modification of the iterative method presented in Aoyama, Sakajo \& Tanaka~\cite{Aoy-Sak-Tan13}. \\
\noindent{\bf Initialization:}\\
Set
\[
z^0_j=\zeta_j, \quad a^0_j=(1-0.5r)\ell_j, \quad b^0_j=r a^0_j,
\]
where $0<r\le1$ is a small positive real number which is the ratio of the lengths of the major and minor axes of the ellipse (see Figure~\ref{f:str} dotted line for $r=0.2$). 

\noindent{\bf Iterations:} \\
For $k=1,2,3,\ldots,$
\begin{itemize}
	\item Use the method presented in Theorem~\ref{thm:cm-rec-slit} to map the preimage domain $G^{k-1}$ to a canonical rectilinear slit domain $\Omega^k$ which is the entire $\zeta$-plane with $m$ slits $L^k_j$, $j=0,1,\ldots,m$, making angles  $\theta_j$ with the positive real axis which are the same as for the given slit domain $\Omega$.
	\item For $j=0,1,\ldots,m$, let $\ell^k_j$ denote the length of the slit $L^k_j$ and let $\zeta^k_j$ denote its center. Then we we define the parameters of the preimage domain $G^k$ as
\begin{eqnarray}
\label{eq:slt-k}
z^{k}_j &=& z^{k-1}_j-(\zeta^{k}_j-\zeta_j), \\
a^{k}_j &=& a^{k-1}_j-(1-0.5r)(\ell^{k}_j -\ell_j), \\
b^{k}_j &=& r a^{k}_j.
\end{eqnarray}
  \item Stop the iteration if 
	\[
	\frac{1}{m+1}\sum_{j=0}^{m}\left(|\zeta^{k}_j -\zeta_j|+|\ell^{k}_j -\ell_j|\right)<\varepsilon \quad{\rm or}\quad k>{\tt Max}
	\]
	where $\varepsilon$ is a given tolerance and ${\tt Max}$ is the maximum number of iterations allowed. In our
numerical calculations we always used $\varepsilon=10^{-14}$ and ${\tt Max}=100$.
\end{itemize}

It is clear that in each iteration of the iterative method, it is required to solve the integral equation with the generalized Neumann kernel~(\ref{e:ie}) and to compute the function $h$ in~(\ref{e:h}) which will be done using the MATLAB function \verb|fbie| as explained in \S\ref{sec:num}. The number of GMRES iterations required for solving the integral equation depends on~$r$. For fixed $r$, the number of GMRES iterations is almost the same for each iteration. Further, as was reported in~\cite{NAS-ETNA15}, the number of GMRES iterations is almost independent of $n$.

\begin{figure}[ht] %
{\quad(a)\hfill\qquad\qquad\qquad\qquad(b)\hfill}

\centerline{
\scalebox{0.5}{\includegraphics{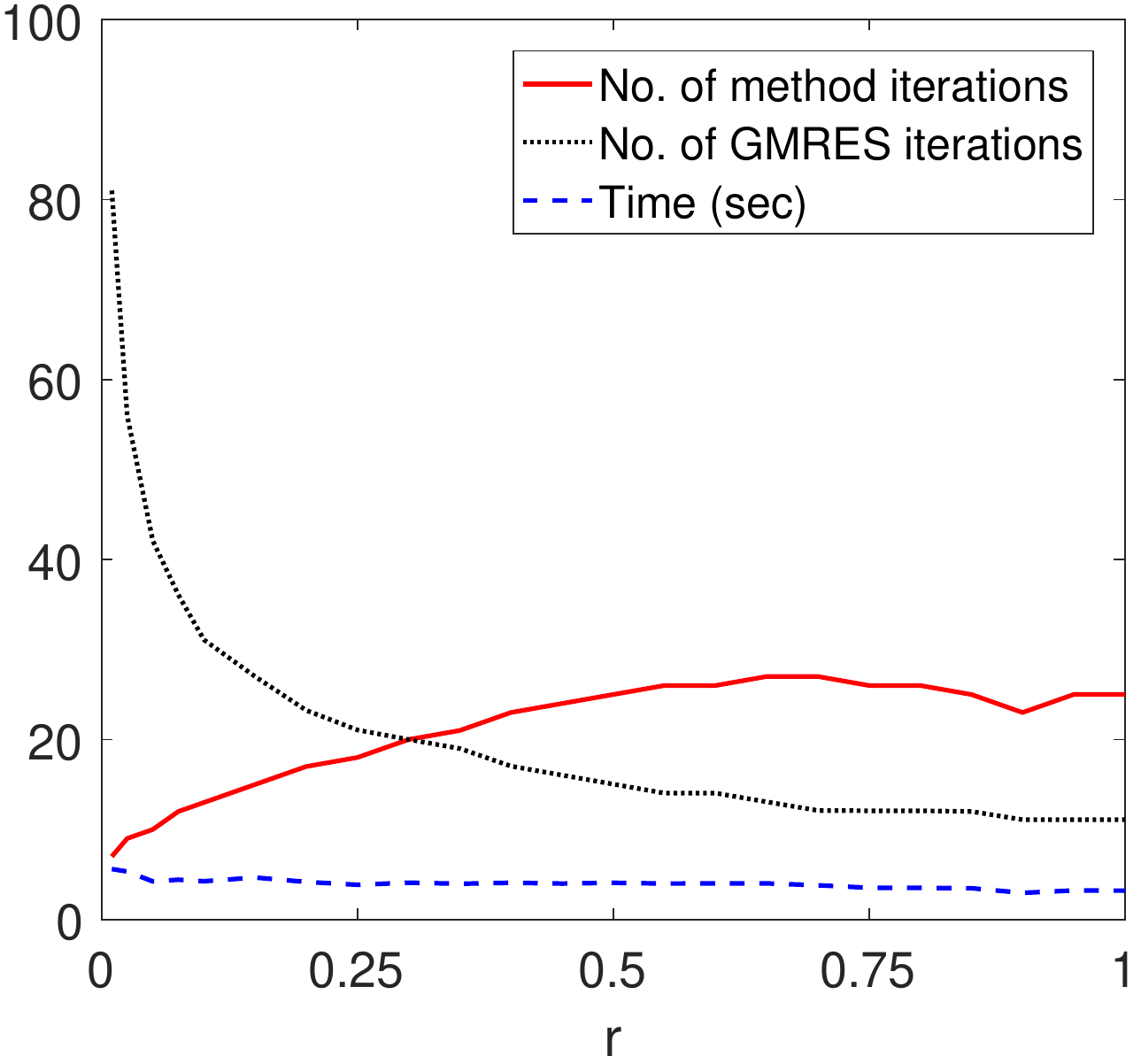}}
\hfill
\scalebox{0.5}{\includegraphics{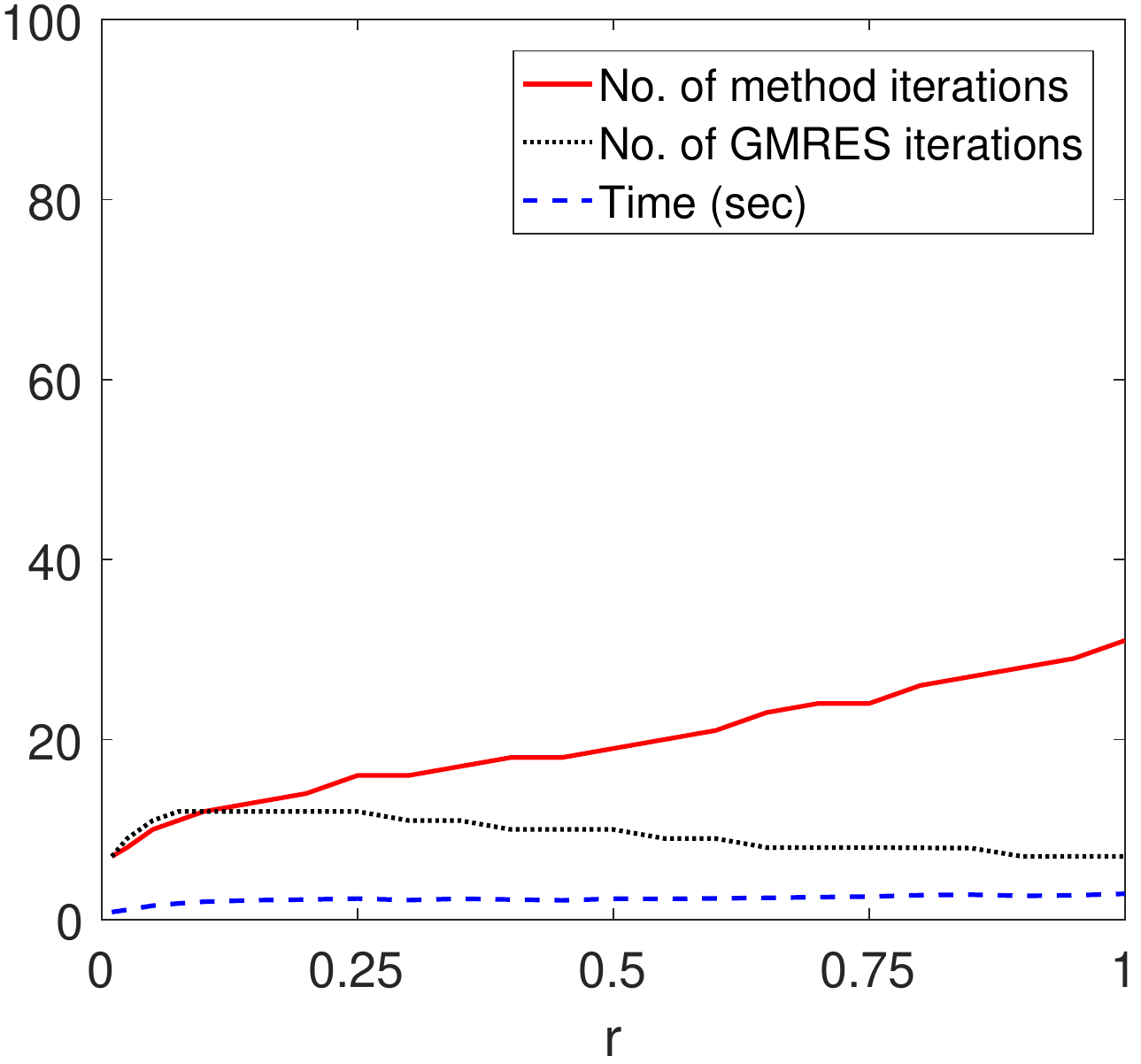}}
}
{\quad(c)\hfill\qquad\qquad\qquad\qquad(d)\hfill}

\centerline{
\scalebox{0.5}{\includegraphics{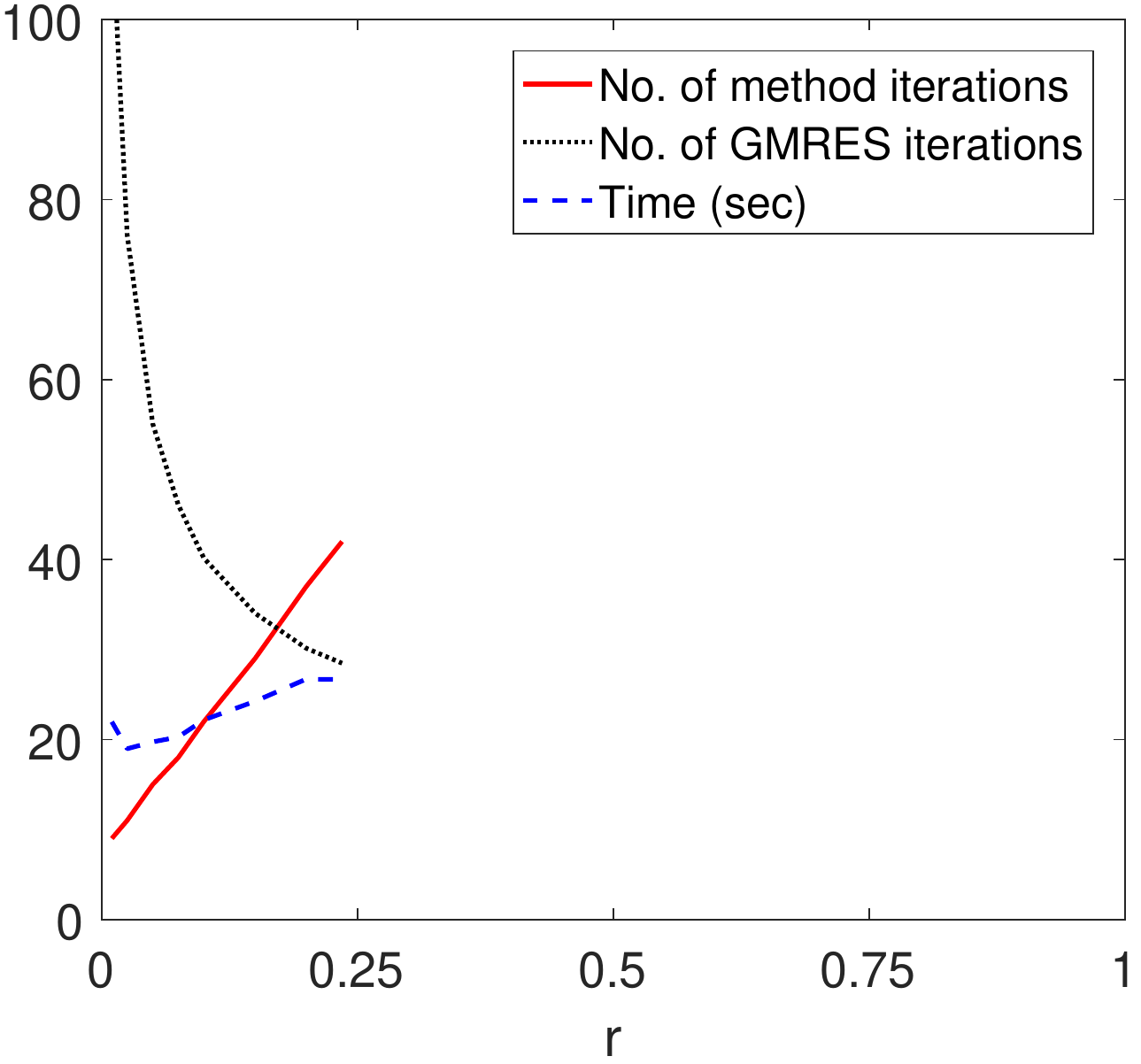}}
\hfill
\scalebox{0.5}{\includegraphics{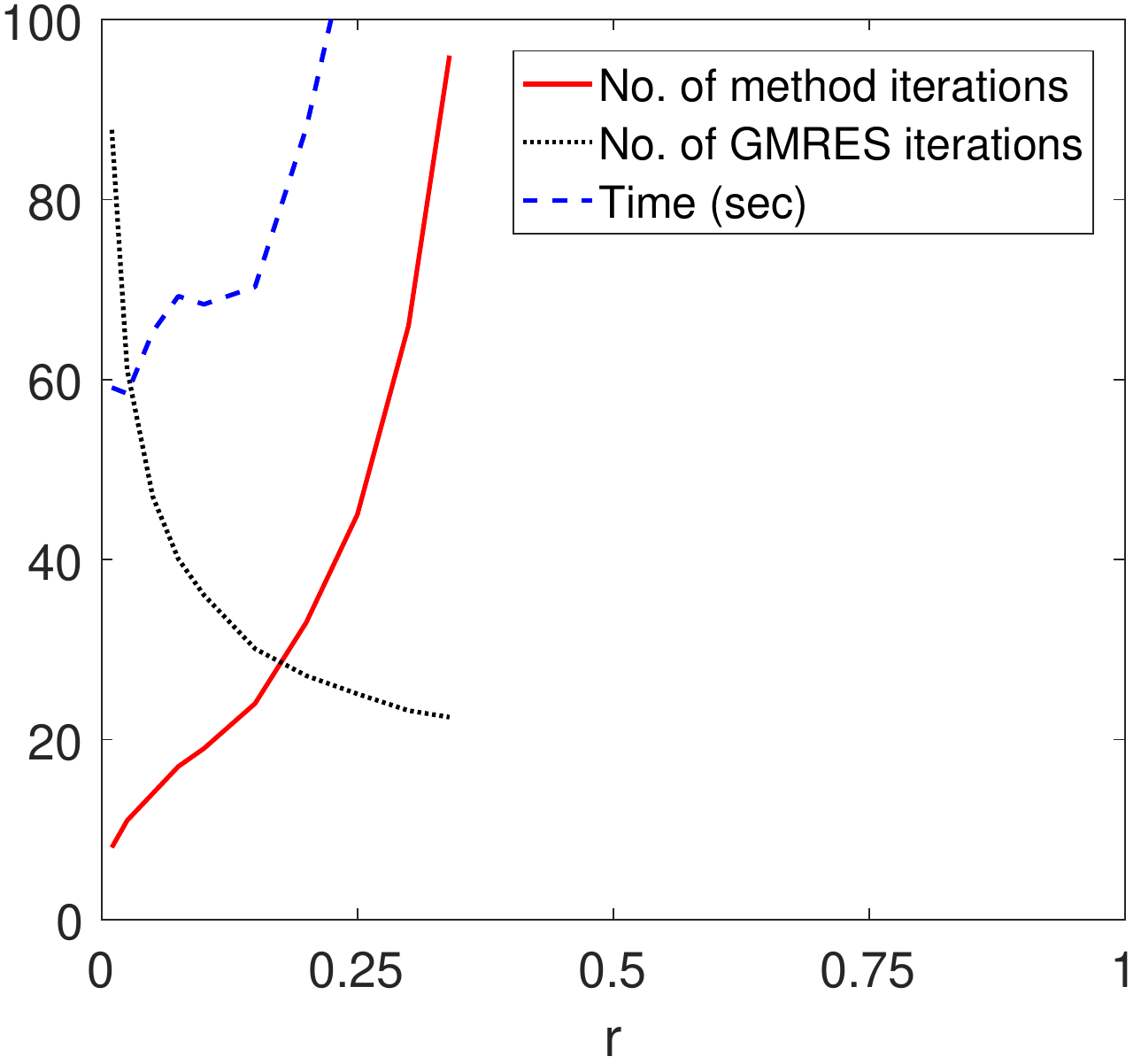}}
}
\caption{The number of iterations required for convergence of the iterative method, the average of the number of GMRES iterations required for solving the integral equation for all iterations, and the total CPU time (in seconds) required to calculate the preimage domain $G$ versus the ratio $r$, for the slit domains shown in Figure~\ref{f:Stirrers-str}. These numerical results are obtained with $n=2^{10}$ nodes per boundary component.} 
\label{f:itr-str}
\end{figure}

\begin{figure}[ht] %
{\quad(a)\hfill\qquad(b)\hfill}

\centerline{
\scalebox{0.6}{\includegraphics{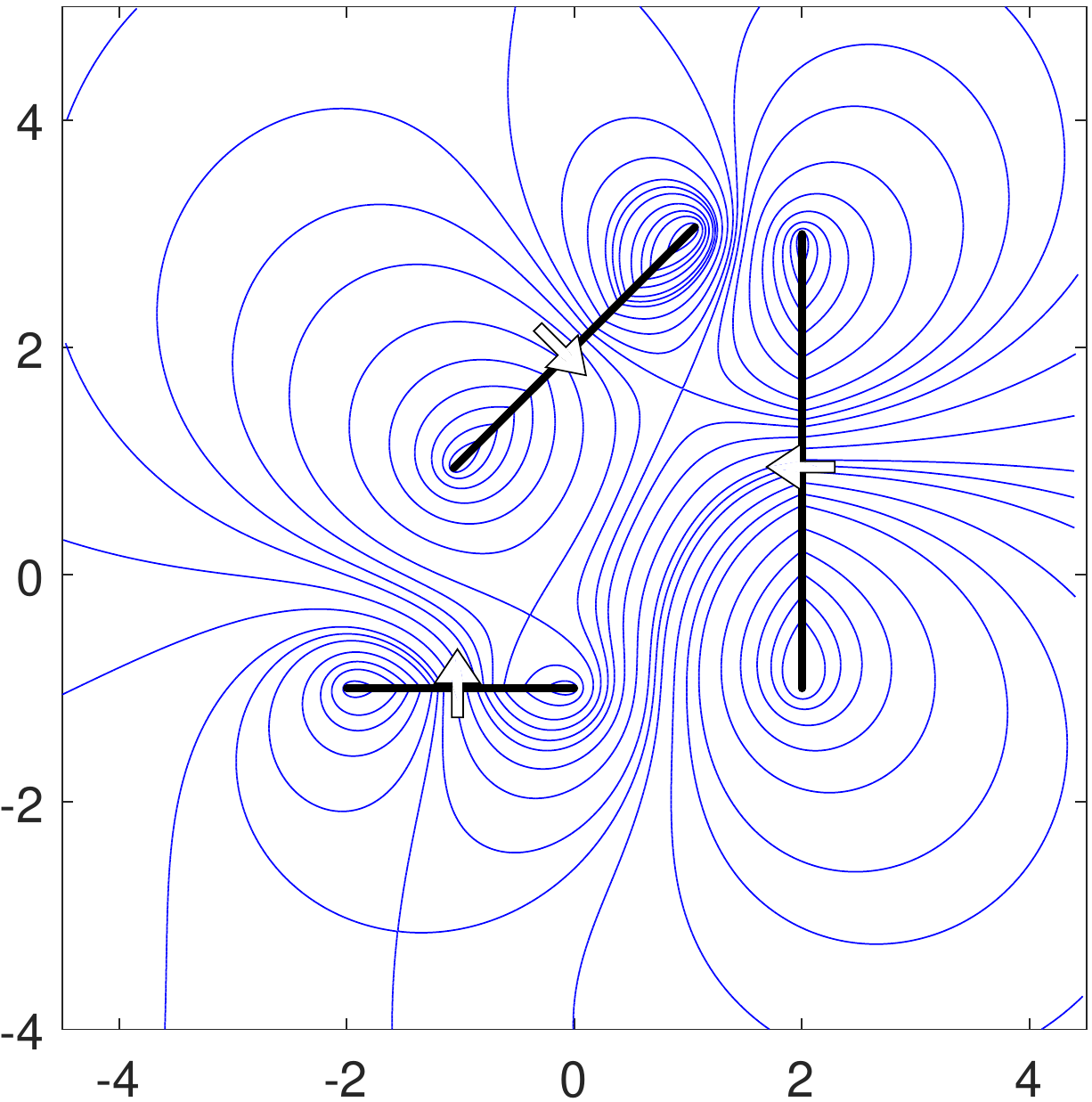}}
\hfill
\scalebox{0.6}{\includegraphics{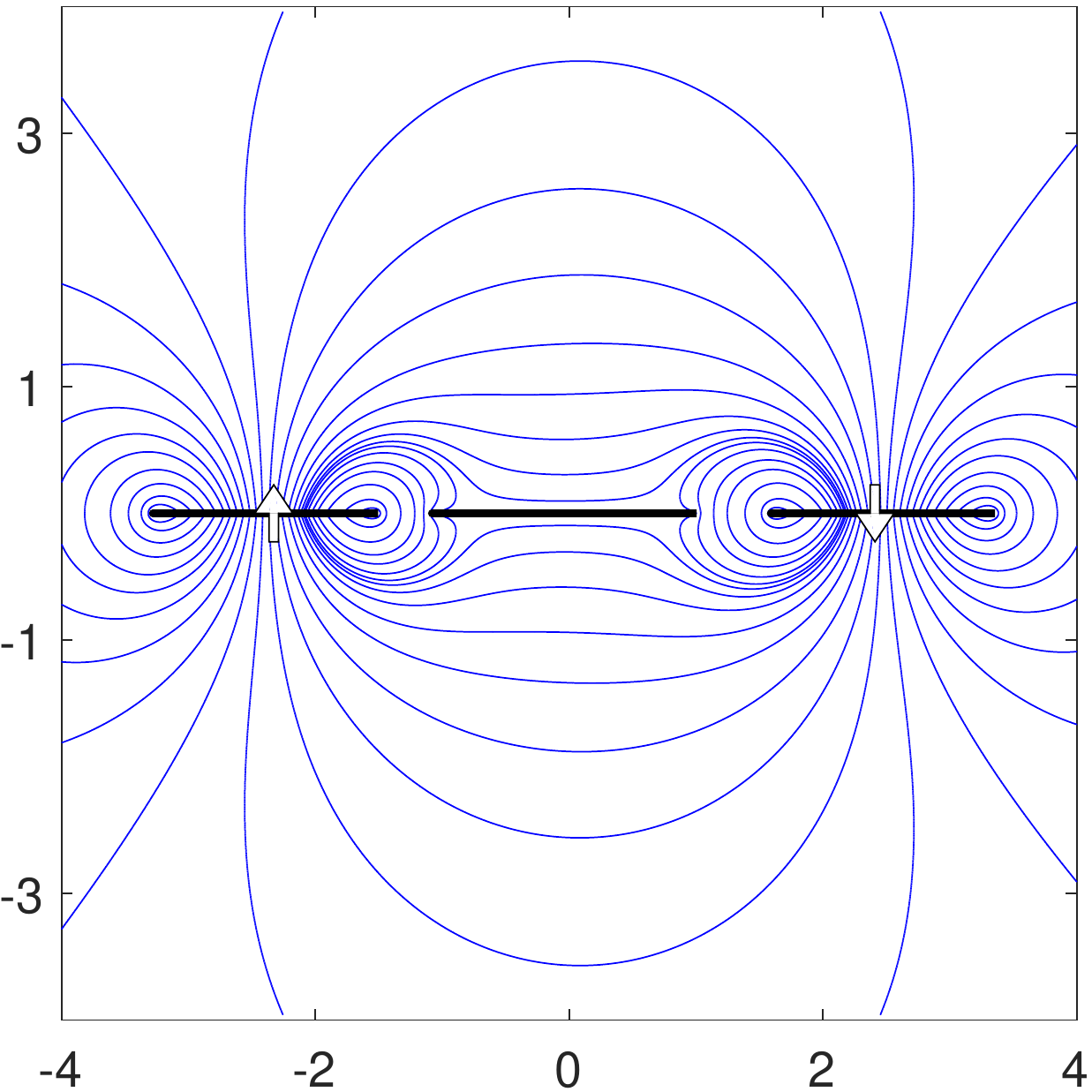}}
}
{\quad(c)\hfill\qquad\qquad(d)\hfill}

\centerline{
\scalebox{0.6}{\includegraphics{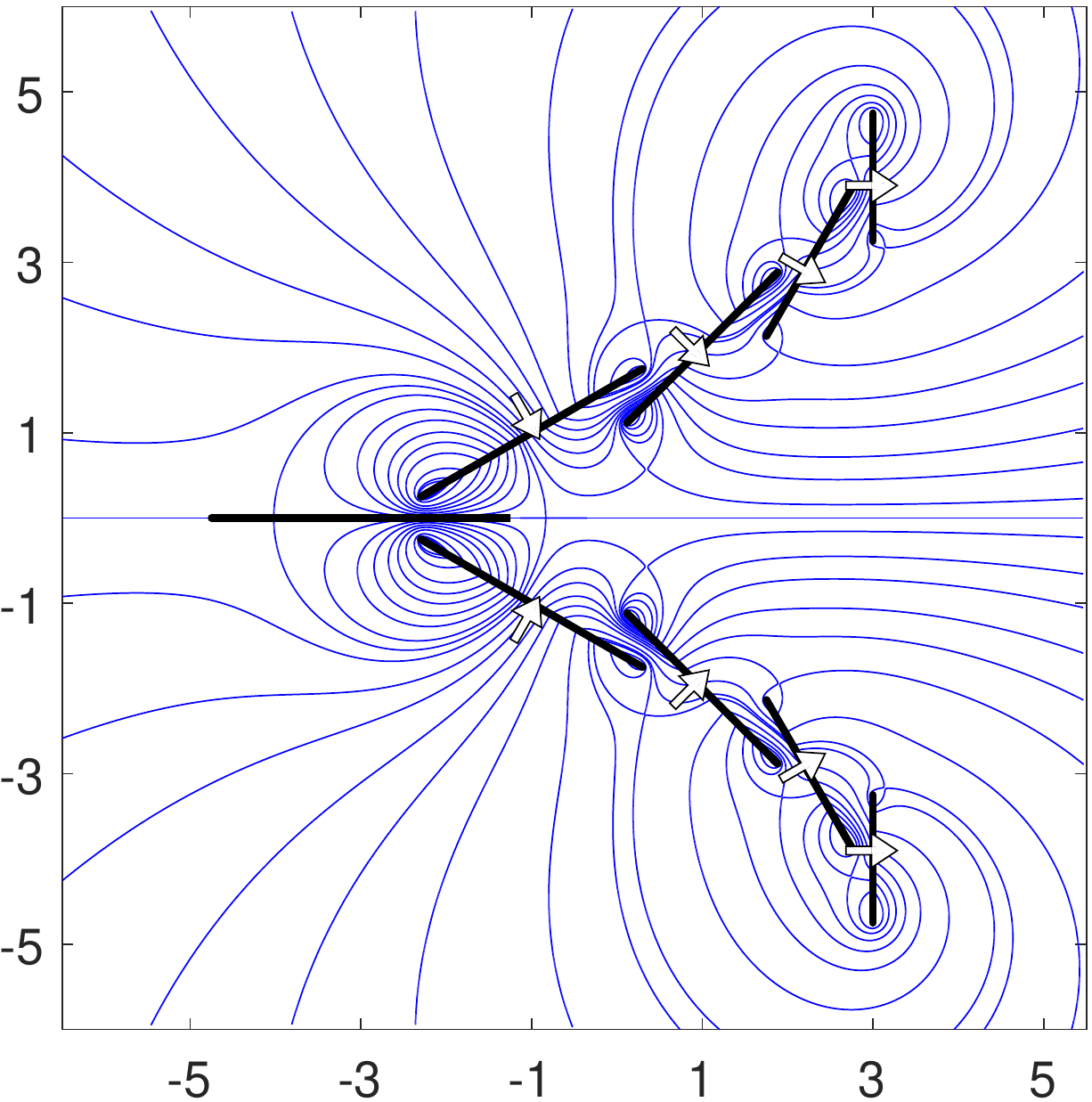}}
\hfill
\scalebox{0.6}{\includegraphics{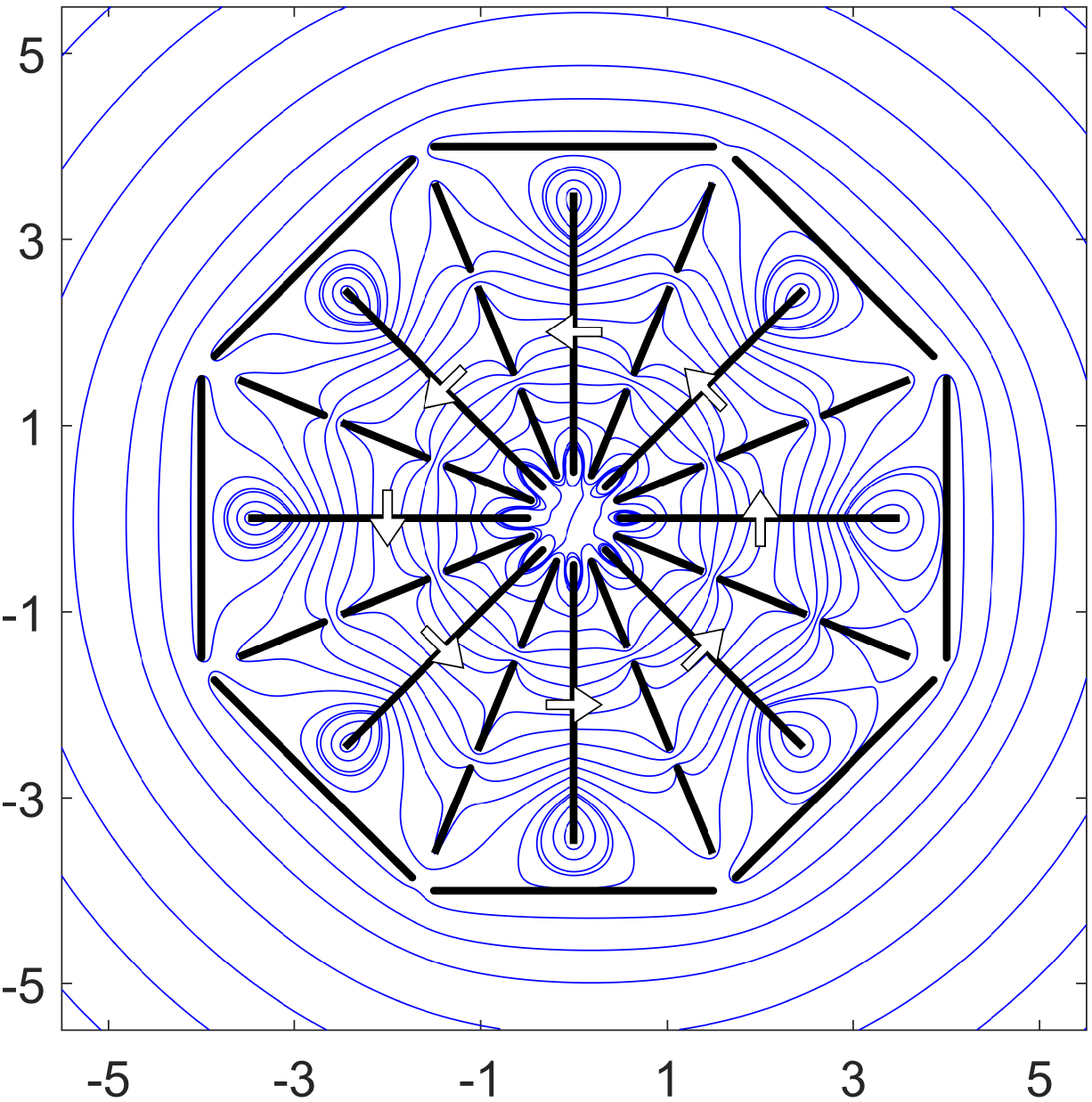}}
}
\caption{Streamlines of the flow generated by rectilinear stirrers in an unbounded domain. Stirrers with arrows have complex velocity of modulus $1$ in the directions indicated by the arrows. Stirrers without arrows are stationary. All slits have zero circulation around them except in the case of  the external slits in (d) which have associated circulation $-1$. The slit domain in (b) has been considered in~\cite[Fig.~6]{cro-str}. Streamlines have been computed with $r=0.2$ and $n=2^{10}$ nodes per boundary component.} 
\label{f:Stirrers-str}
\end{figure}

The above algorithm has been tested for four rectilinear slit domains (see Figure~\ref{f:Stirrers-str}). We assume that $0<r\le1$. For $r=1$, the boundaries $\Gamma_j$ are circles and the preimage domain $G$ is circular.  It turns out that the number of iterations required for convergence of the iterative method increases when $r$ increases. However, the number of GMRES iterations required for solving the integral equation decreases when $r$ increases. 
Hence, the optimal value of $r$ depends on the geometry of the slit domain. Figure~\ref{f:itr-str} shows the number of iterations required for convergence of the iterative method, the average of the number of GMRES iterations required for solving the integral equation for all iterations, and the total CPU time (in seconds) required to calculate the preimage domain $G$ versus the ratio $r$ for the slit domains shown in Figure~\ref{f:Stirrers-str}. Based on the numerical results presented in Figure~\ref{f:itr-str}, when the slits are well separated (Figures~\ref{f:Stirrers-str}(a,b)), the iterative method converges for all $0.01\le r\le1$. However, when the slits are close together, the iterative method converges only for small $r$ (for $0.01\le r\le0.235$ for Figure~\ref{f:Stirrers-str}(c) and for $0.01\le r\le0.34$ for Figure~\ref{f:Stirrers-str}(d)). Thus, we conclude that when the slits are well separated, we can choose $r=1$. But, for slits that are close to each other, we need to choose small $r$. Finally, it is worth mentioning that we made several unsuccessful numerical experiments in the attempt to find an optimal value of $r$ in terms of the minimum distance between the slits. This issue will continue to be investigated in future research.

By obtaining the preimage domain $G$ and the conformal mapping $\Phi$ from $G$ onto $\Omega$, we can calculate the streamlines of the flow generated by the rectilinear stirrers in an unbounded flow as explained above. The streamlines for four rectilinear slit domains obtained with $n=2^{10}$ nodes per boundary component and the ratio $r=0.2$ are shown in Figure~\ref{f:Stirrers-str}.

The previous iterative method provides us the parametrization $\eta(t)$ of the preimage domain $G$ as well as the boundary values $\Phi(\eta(t))$ of the conformal mapping $\Phi$ from $G$ onto $\Omega$. If we are interesting in computing the values of the inverse mapping $\Psi^{-1}$, we need to compute the derivative $\Phi'(\eta(t))$ numerically. Since $\Phi(\eta(t))$ is $2\pi$-periodic, the derivative $\Phi'(\eta(t))$ can be computed accurately by approximating the real and imaginary part of $\Phi(\eta(t))$ by trigonometric interpolating polynomials and then differentiating. 
The inverse mapping $\Phi^{-1}$ has the following Laurent series expansion near $\infty$:
\[
\Phi^{-1}(w)=w+O\left(\frac{1}{w}\right).
\]
Then for $w\in\Omega$, the values of the inverse map $\Phi^{-1}(w)$ can be computed through the Cauchy integral
\begin{equation}\label{eq:rec-Phi-1}
\Phi^{-1}(w)=w+\frac{1}{2\pi\i}\int_{\partial\Omega}\frac{\Phi^{-1}(\zeta)-\zeta}{\zeta-w}d\zeta.
\end{equation}
By using the parametrization $\zeta(t)=\Phi(\eta(t))$ of the boundary $\partial\Omega$, we obtain
\begin{equation}\label{eq:rec-Phi-2}
\Phi^{-1}(w)=w+\frac{1}{2\pi\i}\int_{J}\frac{\eta(t)-\Phi(\eta(t))}{\Phi(\eta(t))-w}\Phi'(\eta(t))\eta'(t)dt.
\end{equation}

\subsection{The upper half-plane with $m$ finite rectilinear slits}
\label{sc:half}

\begin{figure}[ht] %
\centerline{
\scalebox{0.4}{\includegraphics{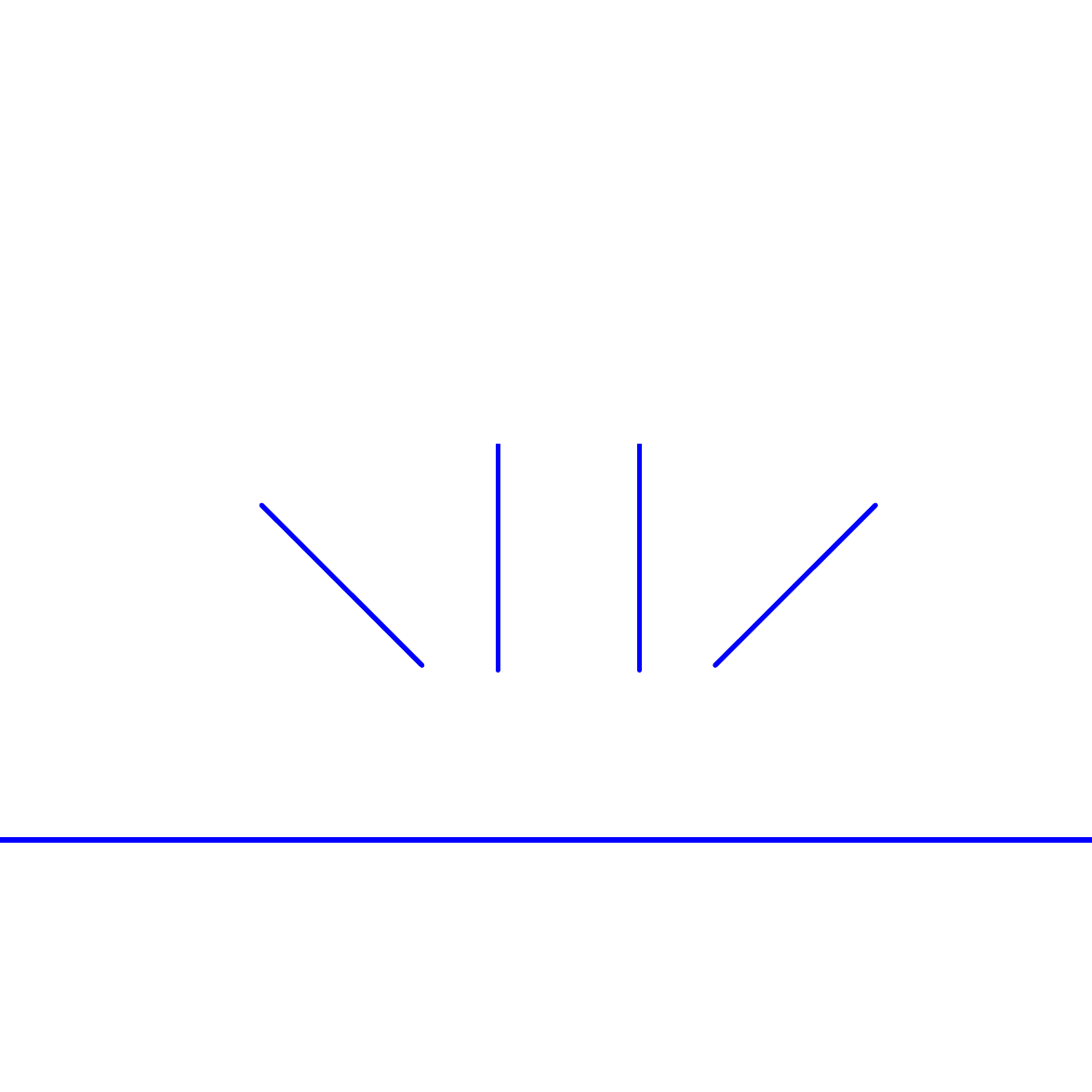}}
\hfill
\scalebox{0.4}{\includegraphics{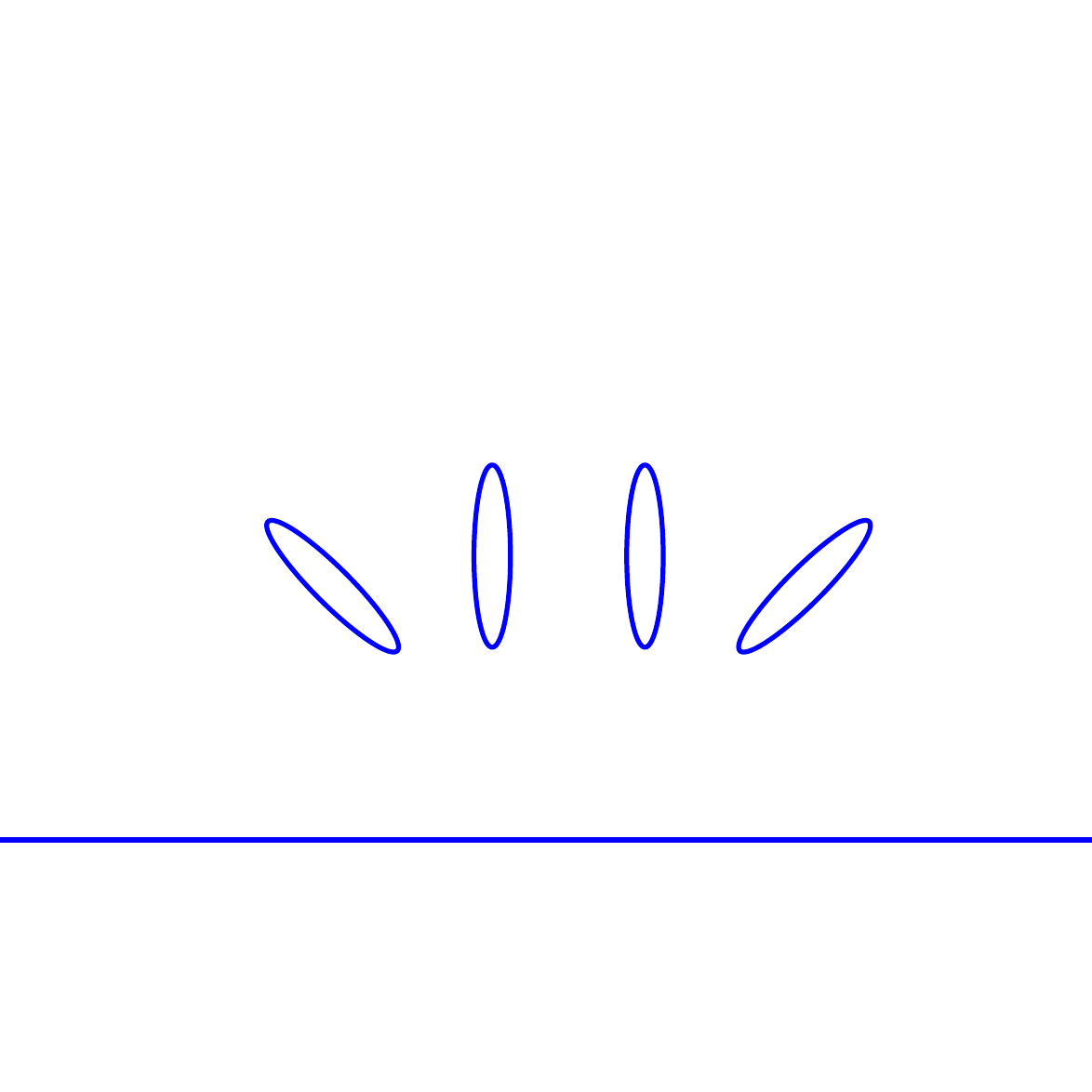}}
\hfill
\scalebox{0.4}{\includegraphics{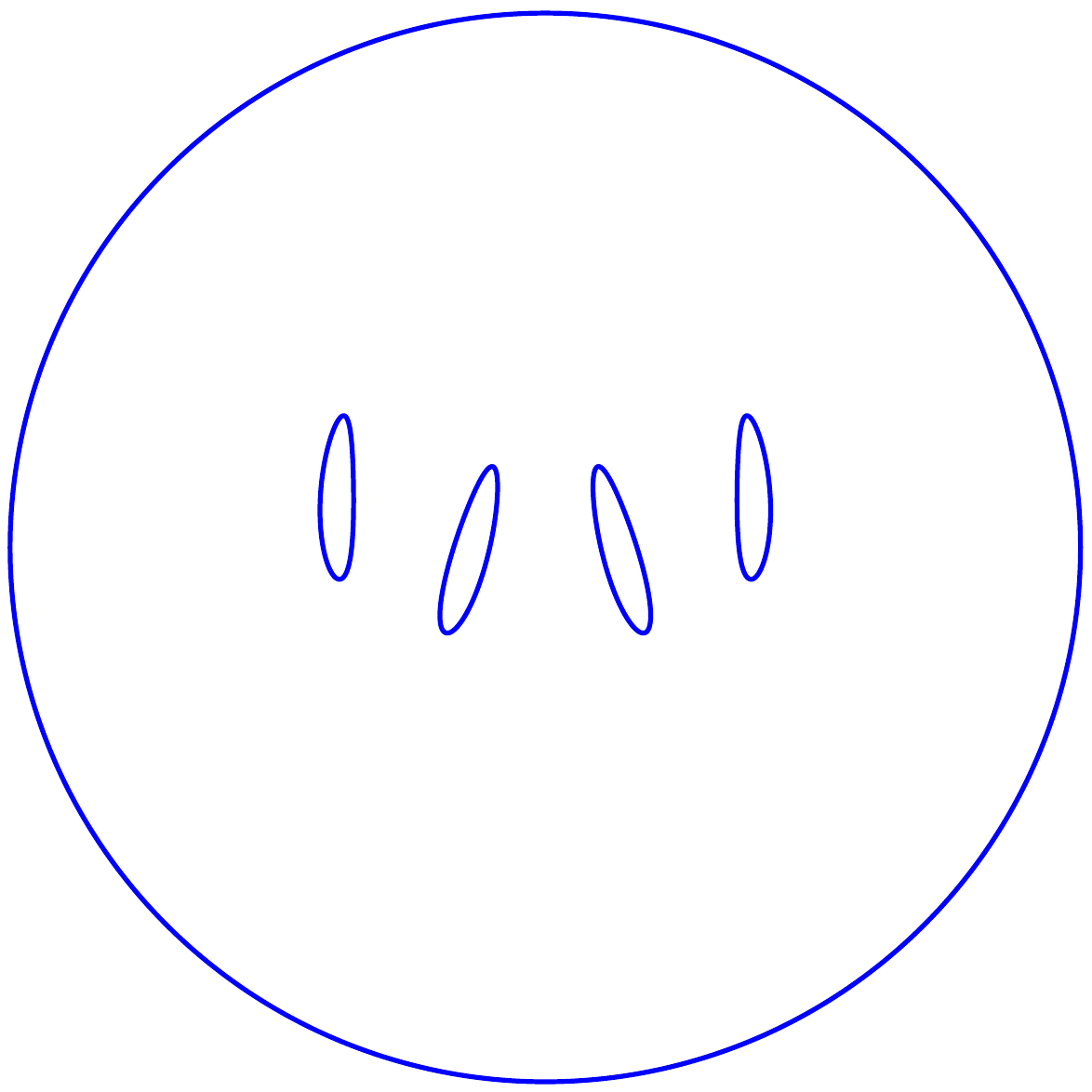}}
}
\caption{The upper half-plane with rectilinear slit domain (left), the initial auxiliary preimage domain $\hat G^0$ (center), and the initial preimage domain $G^0$ (right).} 
\label{f:half}
\end{figure}

This canonical domain consists of the upper half-plane with $m$ rectilinear slits $L_j$, $j=1,2,\ldots,m$ (see Figure~\ref{f:half} (left) for $m=4$). 
For this canonical domain, we will need the following M\"obius transformation  
\[
\xi=\Psi(z)=\i\frac{\i+z}{\i-z}
\]
which maps the unit circle onto the real line and the interior of the unit circle onto the upper-half of the plane with $\Psi(\i)=\infty$ and $\Psi(0)=\i$. Hence, the inverse M\"obius transformation 
\[
z=\Psi^{-1}(\xi)=\i\frac{\xi-\i}{\xi+\i}
\]
maps the real line onto the unit circle and the upper-half of the plane onto the interior of the unit circle. 

To find a preimage domain $G$, we shall consider first an auxiliary preimage domain $\hat G^0$ which is the unbounded multiply connected domain in the upper half-plane and exterior to $m$ ellipses (Figure~\ref{f:half} (center)). 
Thus, the image of the unbounded domain $\hat G^0$ under the mapping $z=\Psi^{-1}(\xi)$ is a bounded domain $G^0$ interior to the unit circle and exterior to $m$ quasi-ellipses (Figure~\ref{f:half} (right)). The domain $G^0$ will be used as an initial approximation of the preimage domain $G$ of the domain $\Omega$ in our numerical calculations. We shall describe an iterative method for computing a sequence of domains $G^0, G^1, G^2,\ldots$ which converges to the preimage domain $G$. In each iteration $k$, it is required to calculate the conformal mapping $\zeta=\Phi(z)$ from $G^{k-1}$ onto a canonical domain $\Omega^k$ which is the upper half-plane with $m$ rectilinear slits $L^k_j$ such that 
\[
\Phi(0)=\i, \quad \Phi(\i)=\infty.
\]
This conformal mapping can be computed as described in the following theorem from~\cite{Nas-JMAA13}.

\begin{theorem}\label{thm:cm-half}
Let $\theta$ be the piecewise constant function defined on $\Gamma$ by $\theta(t)=(0,\theta_1,\ldots,\theta_m)$, the function $A$ be defined by~(\ref{e:A}), and the function $\gamma$ be defined by
\begin{equation}
\gamma(t)=
 \left\{ \begin{array}{l@{\hspace{0.5cm}}l}
0,&t\in J_{0},\\
\Im\left[e^{-\i\theta_j}\Psi(\eta_j(t))\right],&t\in J_j, \quad j=1,2, \ldots,m.
\end{array}
\right.
\end{equation}
Let also $\mu$ be the unique solution of the boundary integral equation~(\ref{e:ie}) and the piecewise constant function $h=(h_0,h_1,\ldots,h_m)$ be given by~(\ref{e:h}). Then the function $f$ with the boundary values
\begin{equation}\label{eq:f-half}
f(\eta(t))=(\gamma(t)+h(t)+\i\mu(t))/A(t)
\end{equation}
is analytic in the bounded domain $G$ and the conformal mapping $\Phi$ is 
given by
\begin{equation}\label{eq:Phi-half}
\Phi(z)=\left[\Psi(z)+zf(z)+\i h_0\right]/(1+h_0), \quad z\in G\cup\Gamma. 
\end{equation}
\end{theorem}

For $j=1,2,\ldots,m$, let $\ell_j$ denote the length of the slit $L_j$, let $\zeta_j$ denote its center, and let $\theta_j$ denote the angle of intersection between the slit and the positive real axis. 
For $k=0,1,2,3,\ldots$, where $k$ denotes the iteration number, we shall assume the preimage domain $G^k$ is the bounded multiply connected domain inside the unit circle parametrized by
\[
\eta^k_0(t)=e^{\i t}, \quad  t\in J_0,
\]
and exterior to $m$ quasi-ellipses $\Gamma_1,\ldots,\Gamma_m$ parametrized by
\[
\eta^k_j(t)=\Psi^{-1}\left(z^k_j+0.5e^{\i\theta^k_j}(a^k_j\cos t-\i b^k_j\sin t)\right), \quad t\in J_j, \quad j=1,2,\ldots,m.
\]
This means $G^k$ is the image under the conformal mapping $z=\Psi^{-1}(\xi)$ of the unbounded multiply connected domain $\hat G^k$ in the upper-half plane $\Im\xi>0$ and exterior to the ellipses $\hat\Gamma_j$ parametrized for $j=1,2,\ldots,m$ by
\[
\hat\eta^k_j(t)=z^k_j+0.5e^{\i\theta^k_j}(a^k_j\cos t-\i b^k_j\sin t), \quad t\in J_j.
\]
The parameters $z^k_j$, $a^k_j$, and $b^k_j$, $j=1,2,\ldots,m$, of the ellipses will be computed using the following iterative method.\\
\noindent{\bf Initialization:}\\
Set
\[
z^0_j=\zeta_j, \quad a^0_j=(1-0.5r)\ell_j, \quad b^0_j=r a^0_j, \quad j=1,2,\ldots,m,
\]
where $0<r\le1$ is the ratio of the lengths of the major and minor axes of the ellipse (see Figure~\ref{f:half} dotted line for $r=0.1$). \\
\noindent{\bf Iterations:} \\
For $k=1,2,3,\ldots$,
\begin{itemize}
	\item Use the method presented in Theorem~\ref{thm:cm-half} to map the preimage domain $G^{k-1}$ to the canonical domain $\Omega^k$ which is the upper-half plane $\Im\zeta>0$ with $m$ rectilinear slits $L^k_j$, $j=1,2,\ldots,m$, making angles $\theta_j$ with the positive real axis.
	\item For $j=1,2,\ldots,m$, let $\ell^k_j$ denote the length of the slit $L^k_j$ and let $\zeta^k_j$ denote its center. Then we update the parameters of the preimage domain $G^k$ as
\begin{eqnarray}
\label{eq:half-k}
z^{k}_j &=& z^{k-1}_j-(\zeta^{k}_j-\zeta_j), \\
a^{k}_j &=& a^{k-1}_j-(1-0.5r)(\ell^{k}_j -\ell_j), \\
b^{k}_j &=& r a^{k}_j.
\end{eqnarray}
  \item Stop the iteration if 
	\[
	\frac{1}{m}\sum_{j=1}^{m}\left(|\zeta^{k}_j -\zeta_j|+|\ell^{k}_j -\ell_j|\right)<\varepsilon \quad{\rm or}\quad k>{\tt Max}
	\]
	where $\varepsilon$ is a given tolerance and ${\tt Max}$ is the maximum number of iterations allowed. In our
numerical calculations we always used $\varepsilon=10^{-14}$ and ${\tt Max}=100$.
\end{itemize}
 
\begin{figure}[ht] %
{\quad(a)\hfill\qquad(b)\hfill}

\centerline{
\scalebox{0.525}{\includegraphics{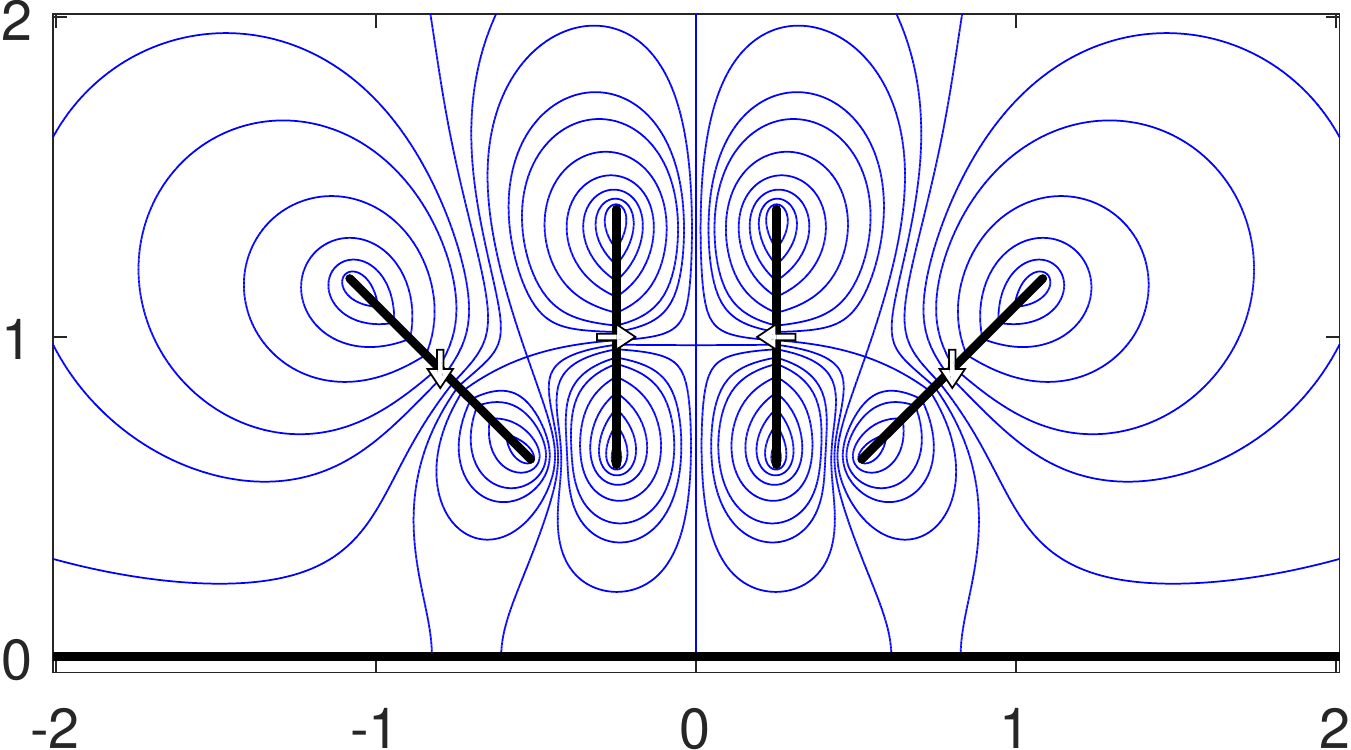}}
\hfill
\scalebox{0.525}{\includegraphics{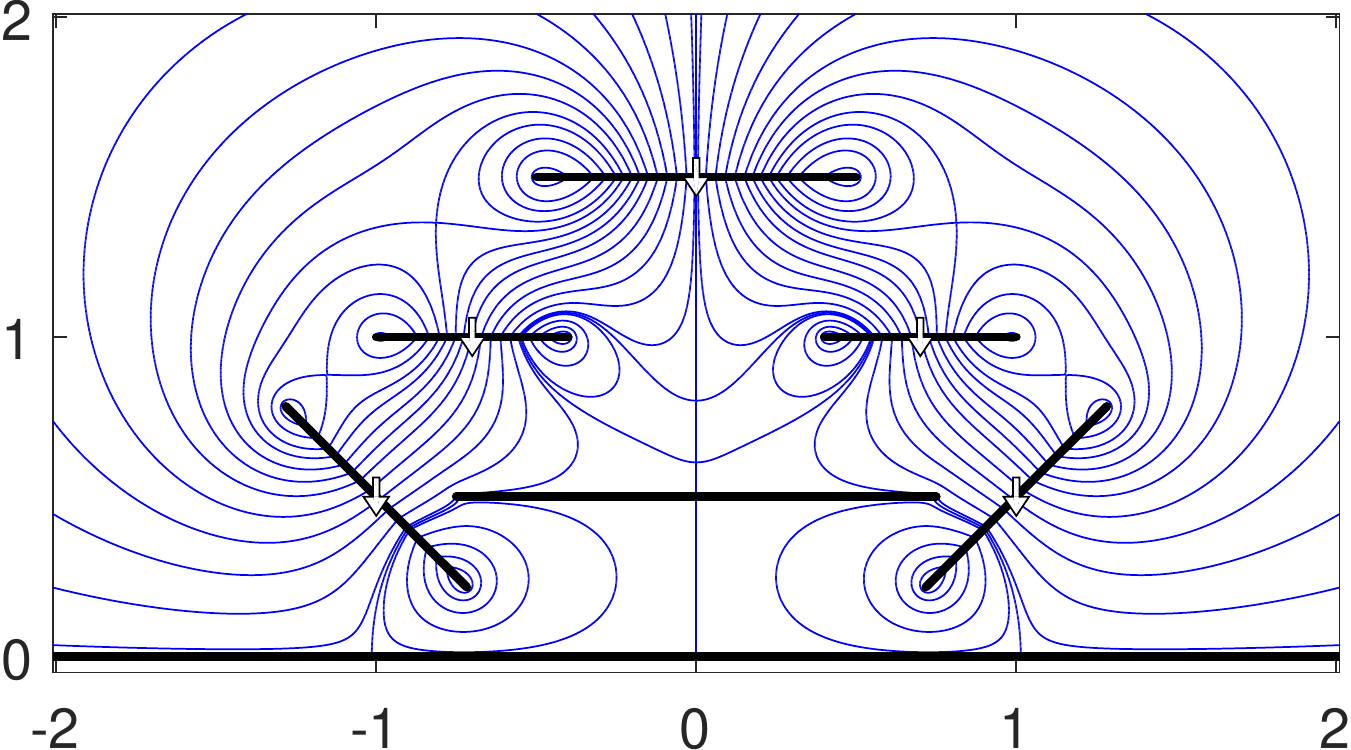}}
}
{\quad(c)\hfill\qquad(d)\hfill}

\centerline{
\scalebox{0.475}{\includegraphics{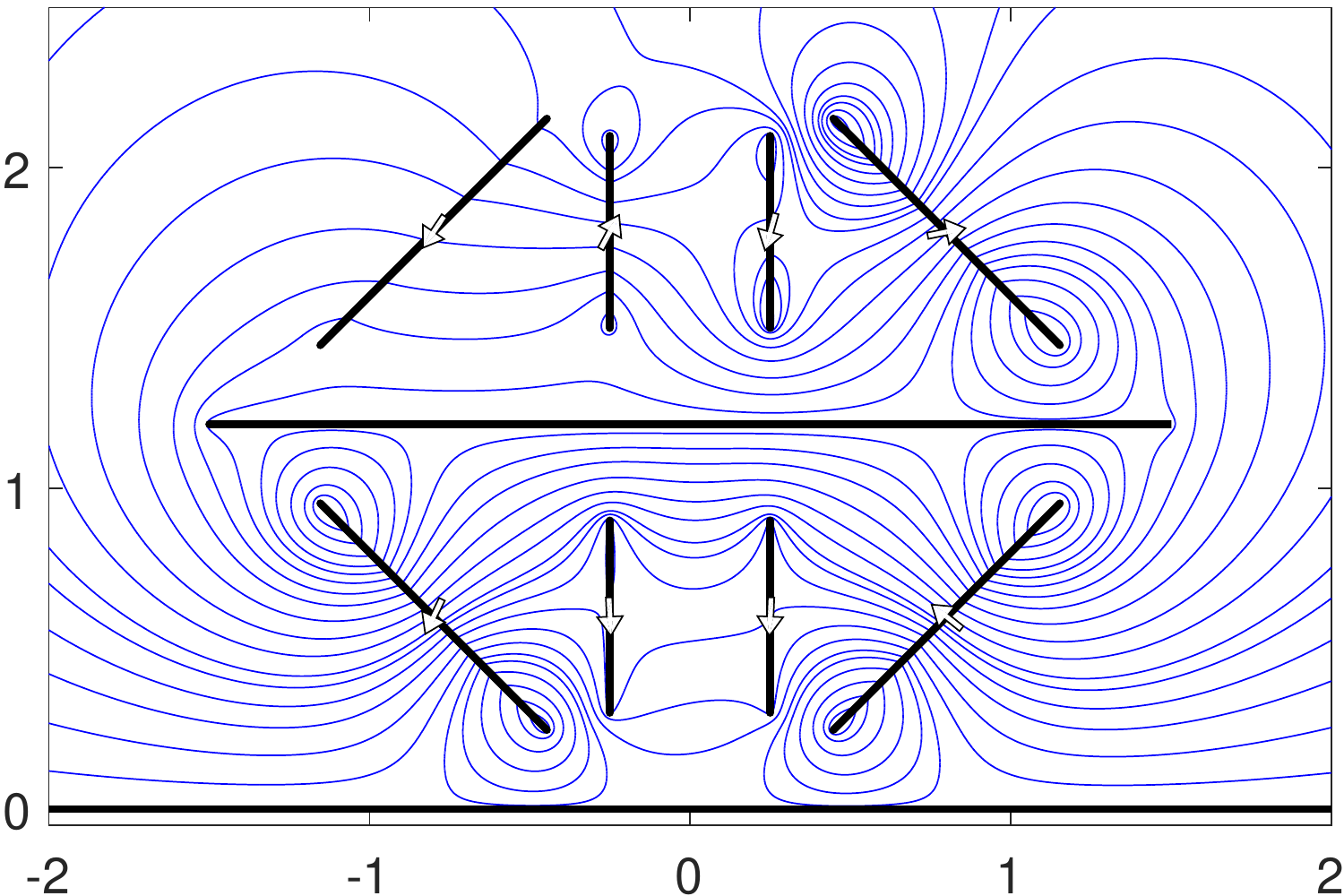}}
\hfill
\scalebox{0.475}{\includegraphics{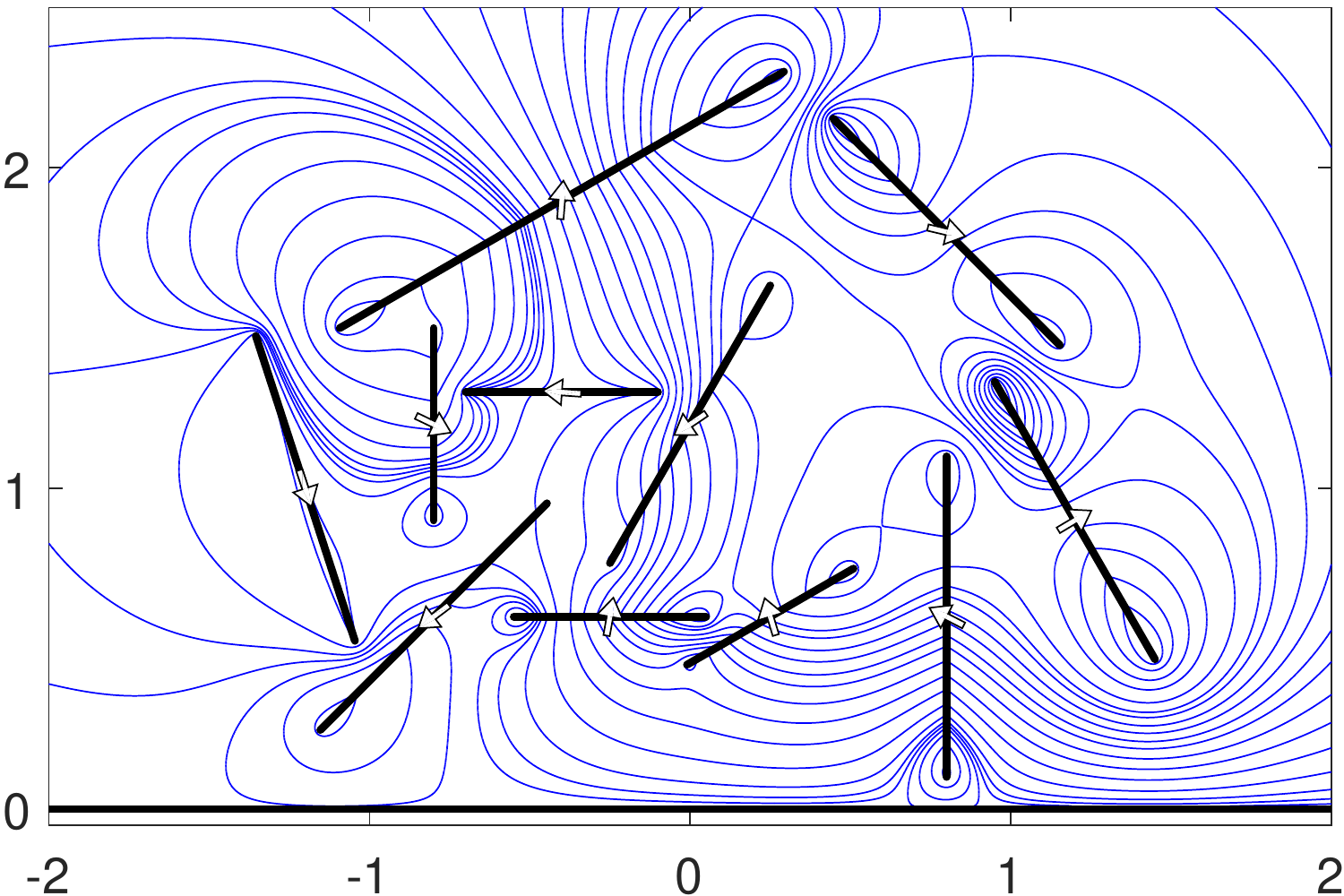}}
}
\caption{Streamlines of the flow generated by rectilinear stirrers in unbounded domains. Stirrers with arrows have complex velocities of modulus $1$ in the directions indicated by the arrows. Stirrers without arrows are stationary. In (a) and (b), all stirrers have zero circulation around them. In (c), all stirrers have zero circulation around them except for the horizontal stirrer which has circulation $-1$ around it. The stirrers in (d) have random circulations between $-1$ and $+1$ around them. These streamlines have been computed with $r=0.1$ and $n=2^{10}$ nodes per boundary component.} 
\label{f:Stirrers-half}
\end{figure}

The algorithm will be tested for four half-plane with rectilinear slit domains (see Figure~\ref{f:Stirrers-half}). By obtaining the preimage domain $G$ and the conformal mapping $\Phi$ from $G$ onto $\Omega$, we calculate the streamlines of the irrotational flow generated by the rectilinear stirrers in an unbounded flow in the above half-plane as explained above. The streamlines obtained with $n=2^{10}$ nodes points per boundary component and the ratio $r=0.1$ are shown in Figure~\ref{f:Stirrers-str}. 

If we are interesting in computing the values of the inverse mapping $\Psi^{-1}(w)$ for $w\in\Omega$, then we can compute these values numerically as long as the values of $\Phi(\eta(t))$ are known.  
Since one of the boundaries of $\Omega$ is unbounded (the whole real line), so instead of computing directly the inverse mapping $\Phi^{-1}(w)$, we shall compute the analytic function $F$ in the domain $\hat\Omega$ defined by
\begin{equation}\label{eq:F-half-1}
F(\xi)=\Phi^{-1}(\Psi(\xi)), \quad \xi\in\hat\Omega,
\end{equation}
where $\hat\Omega$ is the image of the domain $\Omega$ under the M\"obius transform $\Psi^{-1}$ (note that all boundaries of $\hat\Omega$ are bounded). The boundary $\partial\Omega$ of the domain $\Omega$ is parametrized by $\zeta(t)=\phi(\eta(t))$. Hence, the boundary $\partial\hat\Omega$ of $\hat\Omega$ is parametrized by
\begin{equation}\label{eq:h-zeta}
\hat\zeta(t)=\Psi^{-1}(\zeta(t))=\Psi^{-1}(\Phi(\eta(t))), \quad t\in J.
\end{equation}
Then by the Cauchy integral formula, we have
\begin{equation}\label{eq:F-half-2}
F(\xi)=\frac{1}{2\pi\i}\int_{\partial\hat\Omega}\frac{F(\hat\zeta)}{\hat\zeta-\xi}d\hat\zeta.
\end{equation}
Then by using the parametrization $\hat\zeta(t)=\Psi^{-1}(\Phi(\eta(t)))$, $t\in J$, of the boundary $\partial\hat\Omega$, we obtain
\begin{equation}\label{eq:F-half-3}
F(\xi)=\frac{1}{2\pi\i}\int_{J}\frac{F(\Psi^{-1}(\Phi(\eta(t))))}{\zeta(t)-\xi}\,\hat\zeta'(t)dt,
\end{equation}
where the values of $\hat\zeta'(t)$ can be computed numerically as explained at the end of \S\ref{sc:rec}.
By the definition of the function $F$, we have $F(\Psi^{-1}(\Phi(\eta(t))))=\eta(t)$. Hence the function $F$ can be computed for all $\xi\in\hat\Omega$ through
\begin{equation}\label{eq:F-half-4}
F(\xi)=\frac{1}{2\pi\i}\int_{J}\frac{\eta(t)}{\zeta(t)-\xi}\,\hat\zeta'(t)dt.
\end{equation}
Consequently, it follows from~(\ref{eq:F-half-1}) that the inverse mapping $\Phi^{-1}$ can be computed for all $w\in\Omega$ by
\begin{equation}\label{eq:Phi-1}
\Phi^{-1}(w)=F(\Psi^{-1}(w))
=\frac{1}{2\pi\i}\int_{J}\frac{\eta(t)}{\zeta(t)-\Psi^{-1}(w)}\,\hat\zeta'(t)dt.
\end{equation}

\section{Conclusions}

In this paper, we studied the problem of fluid stirrers in planar
domains containing ideal fluid: more specifically, we solved a certain
class of R-H problem to determine the fluid motion driven by
collections of rigid stirrers moving at constant speeds. We have seen
through our presented examples that several stirrers, comprising
various shapes, can be used to generate rather complex flow patterns.
We have shown that we were able to deal with complicated
configurations of fluid stirrers, i.e. highly multiply connected
domains, largely due to the efficacy of our numerical scheme. We
employed a proven fast and accurate boundary integral equation 
with the generalized Neumann kernel method
which has also been successful in
generating numerous solutions to various conformal mapping and
potential theory problems (\cite{Nas-CMFT09,Nas-SIAM09,Nas-JMAA11,Nas-JMAA13,Nas-SIAM13,nas-lap}; see also~\cite{NAS-ETNA15} for a review).

We showed in the particular case of circular stirrers that our results
for the streamlines are in good qualitative agreement with those of
other researchers~\cite{cro-str,fin}. The results presented in this 
paper will thus complement these existing works; they are also expected to be of
particular interest to those wishing to gain qualitative insight into
the fluid mechanics associated with stirring, and to those in industry
designing efficient batch stirring devices for various applications.
To demonstrate the versatility of our numerical scheme, we considered
the ideal fluid flow generated by a collection of arbitrary-shaped
stirrers made-up of piecewise smooth boundary curves, and also by a
high number of stirrers. Stirrers having general shapes were
considered because we were able to proceed simply by providing a
uniform discretization tracing-out their boundary curves (i.e. without
knowledge of a conformal mapping), giving us the freedom to work over
any fluid domain we wish. In the case of stirrers of slit type, we
presented an effective way to still use our numerical scheme by first
proceeding through conformally equivalent elliptical or
quasi-elliptical domains, an approach first introduced 
in~\cite{Aoy-Sak-Tan13}. We note
that there are analytical formulae, expressed in terms of the
Schottky-Klein prime function, for the conformal mappings to Koebe's~\cite{Koe18}
first category
of canonical multiply connected slit domains (Crowdy \&
Marshall~\cite{cro-mar}). The slit domains considered by us (such as those in
Figures~\ref{f:Stirrers-str} and~\ref{f:Stirrers-half}) 
were very arbitrary and to the best of our knowledge, no such 
explicit conformal mapping formulae exist to these slit domains.

Crowdy~\cite{cro-str} found explicit formulae for the ideal flow due to any
finite number of arbitrarily-shaped fluid stirrers, and he presented
several examples of the induced flow field. He did not undertake
computations for stirrer domains of the same variety as we have
considered in this paper due to computational restrictions related to
the Schottky-Klein prime function defined over highly multiply
connected circular domains at the time; however, new effective software is now
available to compute this special function if it is required in problems where the domains are highly multiply connected~\cite{skpf}. The formulae in~\cite{cro-str} also require knowledge of conformal
maps from multiply connected circular domains to complicated target domains, and these are not always easy or even possible to establish (e.g. those
comprised of boundary curves of differing shapes, like those we
presented in Figures~\ref{f:bou-unbou-15}, \ref{f:Squares-44}, \ref{f:Stirrers-str} and~\ref{f:Stirrers-half}). 
There is no doubt that having the
explicit formulae of~\cite{cro-str} for the problem of fluid stirrers is extremely
valuable, but what we have offered in this paper is an effective
alternative approach which can be used to generate accurate numerical
solutions to this problem. It has the particular advantage of being computationally inexpensive and can be used with minimal geometrical restrictions on the target fluid domain, in
addition to being especially useful when dealing with flow domains with many fluid stirrers.

\section*{Acknowledgements}

MMSN and CCG both acknowledge financial support from Qatar University
grant QUUG-CAS-DMSP-15$\backslash$16-27. CCG acknowledges support from Australian
Research Council Discovery Project DP140100933; he is also grateful
for the hospitality of the Department of Mathematics, Statistics \&
Physics at Qatar University where this work was completed.


\end{document}